\RequirePackage{fix-cm}
\documentclass[smallextended, envcountsame]{svjour3}
\smartqed
\journalname{}

\usepackage[unicode = false,
    pdftoolbar = true,
    colorlinks = true,
    linkcolor = blue,
    citecolor = blue,
    filecolor = black,
    urlcolor = blue,
    breaklinks = true]{hyperref}

\usepackage[utf8]{inputenc}
\usepackage{listings}  % For Code 
\usepackage{enumitem}  
\usepackage{amsmath}   % For Math, duh
\usepackage{amssymb}   % For mathsfrs
\usepackage{xcolor}    % For Code Color
\usepackage{graphicx}  
\usepackage{verbatim}
\usepackage{geometry}
\usepackage{mathtools}
\usepackage{float}
\usepackage{bm}
\usepackage[ruled]{algorithm2e}
\usepackage{subcaption}
\usepackage{placeins}
\usepackage{listings}
\usepackage[justification=centering]{caption}
\usepackage{multirow}
\usepackage{orcidlink}
\usepackage{bbding}
\usepackage{tikz}
\usetikzlibrary{arrows.meta, positioning, shapes.geometric, shapes.misc, fit, calc}

\definecolor{gab}{HTML}{50c878}
\definecolor{pat}{HTML}{bf68ff}
\definecolor{sham}{HTML}{7faaff}

\allowdisplaybreaks[4]

\newcommand{\myalg}{\mathrm{CLARSTA}}
\newtheorem{assump}{Assumption}
\newtheorem{coroll}{Corollary}
\newtheorem{defini}{Definition}

\newtheorem{lemm}{Lemma}
\newtheorem{rema}{Remark}
\newtheorem{theo}{Theorem}
\newtheorem{prop}{Proposition}

\DeclareMathOperator*{\argmin}{argmin}

\DeclareMathAlphabet{\mymathbb}{U}{BOONDOX-ds}{m}{n}

\newcommand{\appdiam}{\overline{\mathrm{diam}}}

\makeatletter
\newcommand*{\rom}[1]{\expandafter\@slowromancap\romannumeral #1@}
\makeatother

\SetKwIF{If}{ElseIf}{Else}{if}{}{else if}{else}{end if}%
\SetKwRepeat{Do}{do}{while}
\allowdisplaybreaks[4]
\SetArgSty{textnormal}

\title{%Random subspace trust-region algorithm for general convex-constrained derivative-free optimization
CLARSTA: A random subspace trust-region algorithm for convex-constrained derivative-free optimization \thanks{No part of this manuscript has been published or submitted elsewhere, and it will not be submitted elsewhere while the manuscript is under review at Mathematical Programming.}}

\author{Yiwen Chen\,\orcidlink{0000-0002-8720-932X}\,\Envelope\and Warren Hare\,\orcidlink{0000-0002-4240-3903}\and Amy Wiebe\,\orcidlink{0000-0002-0804-6252}}
\institute{Department of Mathematics, University of British Columbia, Kelowna, British Columbia, V1V 1V7, Canada.\\
This research is partially funded by the Natural Sciences and Engineering Research Council of Canada (cette recherche est partiellement financ\'ee par le Conseil de recherches en sciences naturelles et en g\'enie du Canada), Discover Grants \#2023-03555 and \#2024-04643.\\
\email{yiwchen@student.ubc.ca, warren.hare@ubc.ca, amy.wiebe@ubc.ca}}
\date{\today}

\begin{document}

\maketitle

\begin{abstract}
    %Model-based derivative-free optimization (DFO) methods are one major class of DFO methods widely used in practice.  However, these methods are known to struggle in high dimensions.  Recent research tackles this issue by searching for decreases in low-dimensional subspaces sampled randomly in each iteration.  However, existing research on random subspace model-based DFO methods mainly focuses on unconstrained DFO problems.  The theory and algorithms for high-dimensional constrained DFO problems are yet to be developed.

    This paper proposes a random subspace trust-region algorithm for general convex-constrained derivative-free optimization (DFO) problems.  Similar to previous random subspace DFO methods, the convergence of our algorithm requires a certain accuracy of models and a certain quality of subspaces.  For model accuracy, we define a new class of models that is only required to provide reasonable accuracy on the projection of the constraint set onto the subspace.  We provide a new geometry measure to make these models easy to analyze, construct, and manage.  For subspace quality, we use a normalized random orthogonal projection matrix to provide a method to sample subspaces that preserve the first-order criticality measure by a certain fraction with a probability lower bound independent of the ambient dimension.  Based on all these new theoretical results, we present an almost-sure liminf convergence result and an $\mathcal{O}(\epsilon^{-2})$ expected iteration complexity analysis of our algorithm.  Numerical experiments on problems with dimensions up to 10000 demonstrate the reliable performance of our algorithm in high dimensions.
\end{abstract}
\medskip

\medskip


\section{Introduction}
    Derivative-free optimization (DFO) methods are designed for optimization problems where the derivatives of the objective and/or constraint functions are unavailable or expensive to compute~\cite{audet2017derivative,conn2009introduction}.  In general, DFO methods can be categorized into direct-search methods, which search for decreases along a finite number of chosen directions, and model-based methods, which approximate the objective/constraint functions by constructing approximation models.  Two surveys of DFO methods can be found in \cite{dzahini2024revisiting,larson2019derivative}.  DFO methods have applications in diverse fields, including energy, materials science, and computational engineering design \cite{alarie2021two}.  In recent years, DFO methods have gained increasing interest as they have found effective use in various problems in machine learning \cite{feurer2019hyperparameter,ghanbari2017black,li2021survey}.

    Despite the increasing attention and growth in DFO methods, one major challenge remains: efficiency in high dimensions.  As derivative information is unavailable, the number of function evaluations these methods require increases with the problem dimension. Therefore, most current DFO methods are designed for problems with no more than a hundred variables.
    
    In recent years, one technique proven promising to handle high-dimensional problems is iteratively searching for decreases in random subspaces. For example, an asynchronous parallel direct-search algorithm proposed in \cite{audet2008parallel} applies a parallel space decomposition technique to break the problem into reduced subproblems.  Several other works consider randomly generating search directions in direct-search methods \cite{gratton2015direct,gratton2019direct,roberts2023direct}.  Random subspaces have also been incorporated into model-based methods.  Model-based trust-region methods based on random subspaces are proposed for unconstrained deterministic problems \cite{cartis2023scalable,cartis2024randomized,chen2024qfully} and stochastic problems~\cite{dzahini2024stochastic}, building on related ideas from derivative-based optimization \cite{cartis2022randomiseda,cartis2022randomisedb,shao2021random}.  Based on the same algorithm framework, a randomized sketching method is analyzed in \cite{menickelly2023avoiding}.
    
    Another challenge faced by DFO methods is handling constraints.  Direct search algorithms with convergence analysis have been studied for bound-constrained \cite{lewis1999pattern}, linearly constrained~\cite{lewis2000pattern}, and general nonlinear constrained problems \cite{lewis2002globally}.  A notable generalization of pattern search algorithms, called the mesh adaptive direct search algorithms, is proposed and analyzed in \cite{abramson2006convergence,audet2006mesh}.  Compared to direct search methods, less research studies model-based methods for constrained problems, many of which are based on the model-based trust-region framework \cite{conejo2013global,gratton2011active,powell2015fast}.  However, the convergence of most of these methods requires the models to be fully linear, which, as pointed out in \cite[p.362]{larson2019derivative} and exemplified in \cite[p.2565]{hough2022model}, may be impossible to satisfy in the constrained case.  To address this issue, the authors of \cite{hough2022model} generalize the notions of fully linear models \cite{conn2009global} and $\Lambda$-poisedness \cite{conn2008geometry} to arbitrary convex sets.  A model-based trust-region method based on linear interpolation models is proposed in \cite{hough2022model} for general convex-constrained DFO problems, which is later extended to linear regression and quadratic models \cite{roberts2025model}.
    
    In all model-based DFO methods mentioned above, the ones that focus on high dimensions are only designed for unconstrained problems, while the ones that consider constraints do not include any explicit techniques for handling high dimensions.  In fact, to the best of our knowledge, none of the existing research studies model-based DFO methods for general high-dimensional constrained problems.  This paper attempts to close this gap by proposing the Convex-constrained Linear Approximation Random Subspace Trust-region Algorithm ($\myalg$).  In particular, we consider the following constrained optimization problem in high dimensions:
    \begin{equation*}
        \min_{\boldsymbol{x}\in \mathcal{C}} f(\boldsymbol{x}),
    \end{equation*}
    where $f:\mathbb{R}^n\to\mathbb{R}$ and $\mathcal{C}\subseteq\mathbb{R}^n$ is convex, closed, and has a nonempty interior.  Similar to previously mentioned random subspace model-based DFO methods, the convergence of our algorithm requires the models constructed in the subspaces to provide a certain level of accuracy, and the subspaces to have a certain level of quality.

    In terms of model accuracy, inspired by \cite{cartis2023scalable,hough2022model}, we generalize the definition of fully linear models \cite{conn2009global} by taking both the constraint set and the subspace into consideration.  In particular, we define a new class of models that is only required to provide the same level of accuracy as fully linear models on the projection of the constraint set onto the subspace.  We define a new notion to measure the geometry of sample sets, i.e., finite sets of points at which the objective function is evaluated and that are used to construct the models.  This new notion can be viewed as the matrix spectral norm restricted to the projection of the constraint set onto the subspace.  To construct and manage such models, we follow the procedure introduced in \cite{chen2024qfully}, which is based on the condition number of the direction matrix corresponding to the sample set.  As discussed in~\cite{chen2024qfully}, this procedure is easier for management than that based on $\Lambda$-poisedness \cite{conn2008geometry}, as it does not need to construct and maximize Lagrange polynomials over balls.

    The subspaces used in random subspace DFO methods are normally required to preserve the criticality measure by a certain fraction with a fixed probability.  In \cite{cartis2023scalable}, a notion called $\alpha$-well-aligned matrices is given to define such subspaces for unconstrained problems. An $\alpha$-well-aligned matrix defines a subspace that preserves the norm of $\nabla f$ by a fraction corresponding to $\alpha$.  Several techniques of sampling $\alpha$-well-aligned matrices with certain probability lower bounds are provided in~\cite{cartis2022randomiseda,cartis2022randomisedb,dzahini2024stochastic,shao2021random}.  However, no such definitions or results exist for constrained problems.  In this paper, we first define a normalized convex-constrained version of $\alpha$-well-aligned matrices.  In particular, the new notion requires these matrices to preserve a certain fraction of the first-order criticality measure for convex-constrained optimization.  Then, we use moment properties of uniformly sampled orthogonal projection matrices to provide a sampling technique for such matrices.  We show that the resulting matrices satisfy the required preservation property with a probability lower bound independent of the ambient dimension.

    Based on these new results, we prove that the first-order criticality measure converges to zero in the liminf almost surely and present an $\mathcal{O}(\epsilon^{-2})$ bound on the expected number of iterations required to find an iterate with a sufficiently small criticality measure.  Numerical results comparing the performance of $\myalg$ and $\mathrm{COBYLA}$ \cite{Powell1994} demonstrate the efficiency of $\myalg$ in high dimensions.

    We emphasize that most of the theoretical results in this paper are independent of the algorithm.  In particular, while $\myalg$ only uses linear approximation models, the convergence and complexity results are not limited to linear models.  In fact, we intentionally present them in terms of quadratic models (see Section \ref{sec:converg} for details). This naturally covers linear models while highlighting the possibility of incorporating higher-order models in our algorithm.

    The remainder of this paper is structured as follows.  We end this section by introducing the notations and background definitions used throughout this paper.  Section \ref{sec:cqfully} defines the new class of models, called $(\mathcal{C},\boldsymbol{Q})$-fully linear models, and the new geometry measure.  One technique to construct a class of $(\mathcal{C},\boldsymbol{Q})$-fully linear models is presented.  Section \ref{sec:CLARSTA} describes the general framework of $\myalg$.  We show that the models managed by the algorithm are guaranteed to be $(\mathcal{C},\boldsymbol{Q})$-fully linear.  Moreover, a normalized convex-constrained version of $\alpha$-well-aligned matrices is given, and one method to sample such matrices is provided.  Section \ref{sec:converg} presents the convergence and complexity analysis of $\myalg$.  Section \ref{sec:numexps} designs numerical experiments to compare the performance of $\myalg$ and $\mathrm{COBYLA}$ \cite{Powell1994} in both low and high dimensions. Section \ref{sec:conclusion} concludes the work covered by this paper and suggests some future directions.

\subsection{Notations and background definitions}
    Throughout this paper, vectors are written in lowercase boldface and matrices in uppercase boldface (e.g., $\boldsymbol{x}\in\mathbb{R}^n$ and $\boldsymbol{A}\in\mathbb{R}^{n\times z}$ where $n,z\ge 2$). The notation $B(\boldsymbol{x};\Delta)$ denotes the closed ball centered at $\boldsymbol{x}$ with radius $\Delta$, and $\mathrm{proj}_{\mathcal{C}}(\cdot)$ denotes the Euclidean projection function on $\mathcal{C}$.  For a matrix $\boldsymbol{A}\in\mathbb{R}^{n\times z}$, $\boldsymbol{A}=[\boldsymbol{a}_1\cdots \boldsymbol{a}_z]$ is its column representation, and $\boldsymbol{a}_i\in \boldsymbol{A}$ means that $\boldsymbol{a}_i$ is a column of $\boldsymbol{A}$.  We denote the $(i,j)$-th element of $\boldsymbol{A}$ by $A_{ij}$, its column space by $\mathrm{col}(\boldsymbol{A})$, and the orthogonal complement of $\mathrm{col}(\boldsymbol{A})$ by $\mathrm{col}(\boldsymbol{A})^\perp$.  The all-zero and all-one vectors in $\mathbb{R}^{n}$ are denoted by $\boldsymbol{0}_n$ and $\boldsymbol{1}_n$, respectively.  The identity matrix of order $n$ is denoted by $\boldsymbol{I}_n$ and the $m\times n$ all-zero matrix by $\boldsymbol{0}_{m\times n}$.  The minimum and maximum singular values of a matrix are denoted by $\sigma_{\min}(\cdot)$ and $\sigma_{\max}(\cdot)$, respectively.  The symbol $\|\cdot\|$ denotes the Euclidean norm of a vector and the spectral norm of a matrix, i.e., $\|\boldsymbol{v}\|=\sqrt{\sum_{i=1}^nv_i^2}$ for $\boldsymbol{v}=[v_1,\ldots,v_n]^\top\in\mathbb{R}^n$ and $\|\boldsymbol{A}\|=\sigma_{\max}(\boldsymbol{A})$ for $\boldsymbol{A}\in\mathbb{R}^{n\times z}$ where $n,z\ge 2$.  The set of natural numbers $\{0,1,2,\ldots\}$ is denoted by $\mathbb{N}$. For $k\in\mathbb{N}$, $f\in C^{k}$ means that all partial derivatives of $f$ up to order $k$ exist and are continuous on $\mathbb{R}^n$.  Moreover, $f\in C^{k+}$ means that $f\in C^{k}$ and all the $k$-th order partial derivatives of $f$ are Lipschitz continuous on $\mathbb{R}^n$.

    We adopt the definition and notation of the generalized simplex gradient introduced in \cite{hare2020calculus}.
    \begin{defini}(Generalized simplex gradient)
        Let $\boldsymbol{x}^0\in\mathbb{R}^n$ and $\boldsymbol{D}=[\boldsymbol{d}_1\cdots \boldsymbol{d}_z]\in\mathbb{R}^{n\times z}$.  The generalized simplex gradient of $f$ at $\boldsymbol{x}^0$ over $\boldsymbol{D}$, denoted by $\nabla_S f(\boldsymbol{x}^0;\boldsymbol{D})$, is defined by $$\nabla_S f(\boldsymbol{x}^0;\boldsymbol{D}) = \left(\boldsymbol{D}^\top\right)^\dagger\delta_f(\boldsymbol{x}^0;\boldsymbol{D}),$$ where $(\cdot)^\dagger$ gives the \textnormal{Moore–Penrose} pseudoinverse \cite{penrose1955generalized} of a matrix and
        \begin{equation*}
            \delta_f(\boldsymbol{x}^0;\boldsymbol{D})=\begin{bmatrix}
            f(\boldsymbol{x}^0 + \boldsymbol{d}_1) - f(\boldsymbol{x}^0)\\
            f(\boldsymbol{x}^0 + \boldsymbol{d}_2) - f(\boldsymbol{x}^0)\\
            \vdots\\
            f(\boldsymbol{x}^0 + \boldsymbol{d}_z) - f(\boldsymbol{x}^0)
        \end{bmatrix}.
        \end{equation*}
    \end{defini}

\section{$(\mathcal{C},Q)$-fully linear models}\label{sec:cqfully}
    In model-based DFO, the models are required to satisfy a certain level of accuracy for the methods to converge.  One widely-used class of models for unconstrained problems is called fully linear models \cite{conn2009global}.  For convex-constrained problems, a generalization of fully linear models called $\mathcal{C}$-pointwise fully linear models is defined in \cite{hough2022model,roberts2025model}.  In \cite{cartis2023scalable}, a class of $\boldsymbol{Q}$-fully linear models is defined for unconstrained problems where the models are constructed in the subspace corresponding to the column space of matrix $\boldsymbol{Q}$.
    
    In this section, we first introduce the definition of $(\mathcal{C},\boldsymbol{Q})$-fully linear models, which generalizes the definition of fully linear models \cite{conn2009global}, $\mathcal{C}$-pointwise fully linear models \cite{hough2022model,roberts2025model}, and $\boldsymbol{Q}$-fully linear models~\cite{cartis2023scalable}.  In order to construct and manage $(\mathcal{C},\boldsymbol{Q})$-fully linear models, we introduce a new geometry measure, which can be viewed as the matrix spectral norm restricted to the projection of the constraint set onto the subspace.  We prove that a class of $(\mathcal{C},\boldsymbol{Q})$-fully linear models can be constructed through linear interpolation in the subspace. 

    Let $\boldsymbol{Q}\in\mathbb{R}^{n\times p}$, with $p\le n$, be a matrix whose columns are an orthonormal basis of a $p$-dimensional subspace of $\mathbb{R}^n$.  We define a class of $(\mathcal{C},\boldsymbol{Q})$-fully linear models as follows, which provides the same level of accuracy as fully linear models \cite{conn2009global} on a set equivalent to the projection of $\mathcal{C}$ onto the affine subspace defined by $\boldsymbol{Q}$ at the point of interest $\boldsymbol{x}$, defined by
    \begin{equation*}
        \boldsymbol{Q}^\top\left(\mathcal{C}-\boldsymbol{x}\right) = \left\{\boldsymbol{Q}^\top\left(\boldsymbol{z}-\boldsymbol{x}\right):\boldsymbol{z}\in\mathcal{C}\right\}.
    \end{equation*}
    \begin{defini}\label{def:CQfullylinear}
        Given $\overline{\Delta}>0$, $\boldsymbol{x}\in \mathcal{C}$, $f\in C^1$, 
        and matrix $\boldsymbol{Q}\in\mathbb{R}^{n\times p}$ with $p$ orthonormal columns, we say that $\mathcal{M}_{\overline{\Delta}}=\{\widehat{m}_{\Delta}:\mathbb{R}^p\to\mathbb{R}\}_{\Delta\in (0,\overline{\Delta}]}$ is a class of $(\mathcal{C},\boldsymbol{Q})$-fully linear models of $f$ at $\boldsymbol{x}$ parameterized by $\Delta$ if there exist constants $\kappa_{ef}>0$ and $\kappa_{eg}>0$ such that for all $\Delta\in (0,\overline{\Delta}]$ and $\widehat{\boldsymbol{s}}\in \boldsymbol{Q}^\top(\mathcal{C}-\boldsymbol{x})$ with $\|\widehat{\boldsymbol{s}}\|\le\Delta$,
        \begin{align}
            \left|f(\boldsymbol{x}+\boldsymbol{Q}\widehat{\boldsymbol{s}})-\widehat{m}_{\Delta}(\widehat{\boldsymbol{s}})\right|&\le\kappa_{ef}\Delta^2,\nonumber\\
            \max\limits_{\substack{\boldsymbol{d}\in \boldsymbol{Q}^\top(\mathcal{C}-\boldsymbol{x})\\ \|\boldsymbol{d}\|\le 1}}\left|\left(\boldsymbol{Q}^\top\nabla f(\boldsymbol{x}+\boldsymbol{Q}\widehat{\boldsymbol{s}})-\nabla \widehat{m}_{\Delta}(\widehat{\boldsymbol{s}})\right)^\top \boldsymbol{d}\right|&\le\kappa_{eg}\Delta\label{ineq:CQfullylinear_graderr}.
        \end{align}
    \end{defini}
    \begin{rema}
        If $\mathcal{C}=\mathbb{R}^n$, then $\boldsymbol{Q}^\top(\mathcal{C}-\boldsymbol{x})=\mathbb{R}^p$ and Inequality \eqref{ineq:CQfullylinear_graderr} becomes
        \begin{equation*}
             \max\limits_{\substack{\boldsymbol{d}\in\mathbb{R}^p\\ \|\boldsymbol{d}\|\le 1}}\left|\left(\boldsymbol{Q}^\top\nabla f(\boldsymbol{x}+\boldsymbol{Q}\widehat{\boldsymbol{s}})-\nabla \widehat{m}_{\Delta}(\widehat{\boldsymbol{s}})\right)^\top \boldsymbol{d}\right|=\left\|\boldsymbol{Q}^\top\nabla f(\boldsymbol{x}+\boldsymbol{Q}\widehat{\boldsymbol{s}})-\nabla \widehat{m}_{\Delta}(\widehat{\boldsymbol{s}})\right\| \le \kappa_{eg}\Delta,
        \end{equation*}
        so Definition \ref{def:CQfullylinear} aligns with the definition of $\boldsymbol{Q}$-fully linear models \cite{cartis2023scalable}.
    \end{rema}
    \begin{rema}
        If $\mathcal{C}=\mathbb{R}^n$ and $\boldsymbol{Q}=\boldsymbol{I}_n$, then Definition \ref{def:CQfullylinear} aligns with the classic definition of fully linear models \cite{audet2017derivative,conn2009introduction}, since $\boldsymbol{Q}^\top(\mathcal{C}-\boldsymbol{x})=\mathbb{R}^n$ and
        \begin{equation*}
            \max\limits_{\substack{\boldsymbol{d}\in\mathbb{R}^n\\ \|\boldsymbol{d}\|\le 1}}\left|\left(\nabla f(\boldsymbol{x}+\widehat{\boldsymbol{s}})-\nabla \widehat{m}_{\Delta}(\widehat{\boldsymbol{s}})\right)^\top \boldsymbol{d}\right| = \left\|\nabla f(\boldsymbol{x}+\widehat{\boldsymbol{s}})-\nabla \widehat{m}_{\Delta}(\widehat{\boldsymbol{s}})\right\|.
        \end{equation*}
    \end{rema}
    
    Now, we discuss the construction of a class of $(\mathcal{C},\boldsymbol{Q})$-fully linear models.  For simplicity, we suppose that the model is constructed in the $p$-dimensional affine space
    \begin{equation*}
        \mathcal{Y}=\left\{\boldsymbol{x}^0+\boldsymbol{D}\widehat{\boldsymbol{s}}:\widehat{\boldsymbol{s}}\in\mathbb{R}^p\right\},
    \end{equation*}
    where $p\le n$, $\boldsymbol{x}^0\in\mathcal{C}$ and $\boldsymbol{D}=[\boldsymbol{d}_1\cdots \boldsymbol{d}_p]\in\mathbb{R}^{n\times p}$ has full-column rank.  Suppose that $\boldsymbol{D}=\boldsymbol{Q}\boldsymbol{R}$ is the $\boldsymbol{Q}\boldsymbol{R}$-factorization of $\boldsymbol{D}$, where $\boldsymbol{Q}\in\mathbb{R}^{n\times p}$ and $\boldsymbol{R}=[\boldsymbol{r}_1\cdots \boldsymbol{r}_p]\in\mathbb{R}^{p\times p}$.  Let $\appdiam(\boldsymbol{R})=\max_{1\le i\le p}\|\boldsymbol{r}_i\|$.  Define $\widehat{\boldsymbol{R}}=\boldsymbol{R}\slash\appdiam(\boldsymbol{R})$.   Let $\widehat{f}:\mathbb{R}^p\to\mathbb{R}$ be given by $\widehat{f}(\widehat{\boldsymbol{s}})=f(\boldsymbol{x}^0+\boldsymbol{Q}\widehat{\boldsymbol{s}})$. 
    
    We define a linear model $\widehat{m}:\mathbb{R}^p\to\mathbb{R}$ of $\widehat{f}$ by 
    \begin{equation*}
        \widehat{m}(\widehat{\boldsymbol{s}}) = \widehat{f}(\boldsymbol{0}_p) + \nabla_S\widehat{f}(\boldsymbol{0}_p;\boldsymbol{R})^\top\widehat{\boldsymbol{s}}.
    \end{equation*}
    \begin{prop}
        The model $\widehat{m}(\widehat{\boldsymbol{s}})$ interpolates $\widehat{f}$ on $\{\boldsymbol{0}_p\}\cup\{\boldsymbol{r}_i:i=1,\ldots,p\}$.
    \end{prop}
    \begin{proof}
        Let $\widehat{\boldsymbol{s}}=\boldsymbol{0}_p$.  We have $\widehat{m}(\widehat{\boldsymbol{s}})=\widehat{f}(\boldsymbol{0}_p)$.  Furthermore, notice that
        \begin{equation}\label{eq:GSG(R)multiplyR}
             \nabla_S\widehat{f}(\boldsymbol{0}_p;\boldsymbol{R})^\top \boldsymbol{R} = \delta_{\widehat{f}}(\boldsymbol{0}_p;\boldsymbol{R})^\top \boldsymbol{R}^{-1}\boldsymbol{R} = \delta_{\widehat{f}}(\boldsymbol{0}_p;\boldsymbol{R})^\top,
        \end{equation}
        which implies that $\nabla_S\widehat{f}(\boldsymbol{0}_p;\boldsymbol{R})^\top \boldsymbol{r}_i = \widehat{f}(\boldsymbol{r}_i)-\widehat{f}(\boldsymbol{0}_p)$ for all $i=1,\ldots,p$.  Hence, we have $\widehat{m}(\boldsymbol{r}_i)=\widehat{f}(\boldsymbol{r}_i)$ for all $i=1,\ldots,p$.

    $\hfill\qed$
    \end{proof}
    \begin{rema}
        In fact, the model $\widehat{m}(\widehat{\boldsymbol{s}})$ is a determined linear interpolation model of $\widehat{f}$ in $\mathbb{R}^p$, i.e., $\widehat{m}(\widehat{\boldsymbol{s}})$ is the unique linear model that interpolates $\widehat{f}$ on $\{\boldsymbol{0}_p\}\cup\{\boldsymbol{r}_i:i=1,\ldots,p\}$.  For a detailed discussion, we refer the reader to \cite[Chapters 2-6]{conn2009introduction}.
    \end{rema}
    
    We aim to prove that $\widehat{m}(\widehat{\boldsymbol{s}})$ belongs to a class of $(\mathcal{C},\boldsymbol{Q})$-fully linear models at $\boldsymbol{x}^0$.  We start with the following proposition, which shows that any Lipschitz constant of $\nabla f$ is also a Lipschitz constant of $\nabla\widehat{f}$.
    \begin{prop}
        Suppose that $f\in C^{1+}$ with Lipschitz constant $L_{\nabla f}$.  Then $\widehat{f}\in C^{1+}$ with Lipschitz constant $L_{\nabla f}$.
    \end{prop}
    \begin{proof}
        Since $\boldsymbol{Q}\in\mathbb{R}^{n\times p}$ has $p$ orthonormal columns, we have $\|\boldsymbol{Q}\|=\|\boldsymbol{Q}^\top\|=1$ and $\|\boldsymbol{Q}\boldsymbol{d}\|=\|\boldsymbol{d}\|$ for all $\boldsymbol{d}\in\mathbb{R}^p$.  Let $\widehat{\boldsymbol{s}}_1,\widehat{\boldsymbol{s}}_2\in \mathbb{R}^p$.  Then, we have
        \begin{align*}
            \left\|\nabla\widehat{f}(\widehat{\boldsymbol{s}}_2)-\nabla\widehat{f}(\widehat{\boldsymbol{s}}_1)\right\| 
            &= \left\|\boldsymbol{Q}^\top\nabla f(\boldsymbol{x}^0+\boldsymbol{Q}\widehat{\boldsymbol{s}}_2)-\boldsymbol{Q}^\top\nabla f(\boldsymbol{x}^0+\boldsymbol{Q}\widehat{\boldsymbol{s}}_1)\right\|\\
            &\le \left\|\boldsymbol{Q}^\top\right\|\left\|\nabla f(\boldsymbol{x}^0+\boldsymbol{Q}\widehat{\boldsymbol{s}}_2)- \nabla f(\boldsymbol{x}^0+\boldsymbol{Q}\widehat{\boldsymbol{s}}_1)\right\|\\
            &= \left\|\nabla f(\boldsymbol{x}^0+\boldsymbol{Q}\widehat{\boldsymbol{s}}_2)- \nabla f(\boldsymbol{x}^0+\boldsymbol{Q}\widehat{\boldsymbol{s}}_1)\right\|\\
            &\le L_{\nabla f}\left\|\boldsymbol{Q}\left(\widehat{\boldsymbol{s}}_2-\widehat{\boldsymbol{s}}_1\right)\right\|\\
            &= L_{\nabla f}\left\|\widehat{\boldsymbol{s}}_2-\widehat{\boldsymbol{s}}_1\right\|.
        \end{align*}

    $\hfill\qed$
    \end{proof}
    Henceforth, we do not differentiate the Lipschitz constant of $\nabla f$ and $\nabla\widehat{f}$, and denote them both by $L_{\nabla f}$.

    As pointed out in \cite{hough2022model}, in the constrained case, we may improve the error bounds by restricting the terms related to the geometry of sample sets to the constraint set.  Instead of using $\Lambda$-poisedness \cite{conn2008geometry}, which is used in \cite{hough2022model}, we follow the procedure in \cite{chen2024qfully} and measure the geometry of sample sets based on the matrix spectral norm.  This procedure is easier to manage than that based on the $\Lambda$-poisedness \cite{conn2008geometry} by avoiding maximizing Lagrange polynomials (see \cite{chen2024qfully} for details).  We now define the geometry measure we use in this paper, which can be viewed as the matrix spectral norm restricted to $\boldsymbol{Q}^\top(\mathcal{C}-\boldsymbol{x}^0)$.
    \begin{defini}\label{def:geomeasure}
        Suppose that $\boldsymbol{x}^0\in\mathcal{C}$ and $\boldsymbol{D}\in\mathbb{R}^{n\times p}$ has full-column rank.  Let $\boldsymbol{D}=\boldsymbol{Q}\boldsymbol{R}$ be the $\boldsymbol{Q}\boldsymbol{R}$-factorization of $\boldsymbol{D}$ and denote $\boldsymbol{R}=[\boldsymbol{r}_1\cdots \boldsymbol{r}_p]$.  Let $\appdiam(\boldsymbol{R})=\max_{1\le i\le p}\|\boldsymbol{r}_i\|$.  For a matrix $\boldsymbol{M}\in\mathbb{R}^{m\times p}$, define
        \begin{equation*}
            \left\|\boldsymbol{M}\right\|_{\boldsymbol{x}^0,\mathcal{C},\boldsymbol{D}} = \frac{1}{\min\{\appdiam(\boldsymbol{R}),1\}}\max\limits_{\substack{\boldsymbol{w}\in \boldsymbol{Q}^\top(\mathcal{C}-\boldsymbol{x}^0)\\ \|\boldsymbol{w}\|\le \min\{\appdiam(\boldsymbol{R}),1\}}}\left\|\boldsymbol{M}\boldsymbol{w}\right\|.
        \end{equation*}
    \end{defini}
    \begin{rema}
        For any constant $t\ge 0$, we have $\left\|t\boldsymbol{M}\right\|_{\boldsymbol{x}^0,\mathcal{C},\boldsymbol{D}} = t\left\|\boldsymbol{M}\right\|_{\boldsymbol{x}^0,\mathcal{C},\boldsymbol{D}}$.  For all $\boldsymbol{M}\in\mathbb{R}^{m\times p}$, 
        \begin{equation*}
            \left\|\boldsymbol{M}\right\|_{\boldsymbol{x}^0,\mathcal{C},\boldsymbol{D}}\le\frac{1}{\min\{\appdiam(\boldsymbol{R}),1\}}\max\limits_{\|\boldsymbol{w}\|\le \min\{\appdiam(\boldsymbol{R}),1\}}\left\|\boldsymbol{M}\boldsymbol{w}\right\|=\max\limits_{\|\boldsymbol{w}\|\le 1}\left\|\boldsymbol{M}\boldsymbol{w}\right\|=\left\|\boldsymbol{M}\right\|.
        \end{equation*}  In particular, if $\mathcal{C}=\mathbb{R}^n$, then 
        \begin{equation*}
            \left\|\boldsymbol{M}\right\|_{\boldsymbol{x}^0,\mathcal{C},\boldsymbol{D}}=\frac{1}{\min\{\appdiam(\boldsymbol{R}),1\}}\max\limits_{\|\boldsymbol{w}\|\le \min\{\appdiam(\boldsymbol{R}),1\}}\left\|\boldsymbol{M}\boldsymbol{w}\right\|=\max\limits_{\|\boldsymbol{w}\|\le 1}\left\|\boldsymbol{M}\boldsymbol{w}\right\|=\left\|\boldsymbol{M}\right\|.
        \end{equation*}
    \end{rema}
    \begin{rema}
        We do not need to compute any $\|\boldsymbol{M}\|_{\boldsymbol{x}^0,\mathcal{C},\boldsymbol{D}}$ in the algorithm.  We only need a sample set management procedure that maintains a small value of it (see the error bounds in Subsection~\ref{subsec:modelacc} for details).  In this paper, we use the same sample set management procedure as in \cite{chen2024qfully}, which maintains a small value of $\|\boldsymbol{M}\|$ and hence $\|\boldsymbol{M}\|_{\boldsymbol{x}^0,\mathcal{C},\boldsymbol{D}}$.
    \end{rema}

    Notice that both vectors $\boldsymbol{\widehat{s}}$ and $\boldsymbol{d}$ in Definition \ref{def:CQfullylinear} are in the set $\boldsymbol{Q}^\top(\mathcal{C}-\boldsymbol{x})$.  The next lemma is crucial to develop error bounds for $\widehat{m}(\widehat{s})$, which links the coordinates of a vector in $\boldsymbol{Q}^\top(\mathcal{C}-\boldsymbol{x}^0)$ with length bounded by $\min\{\appdiam(\boldsymbol{R}),1\}$ to $\|\boldsymbol{R}^{-1}\|_{\boldsymbol{x}^0,\mathcal{C},\boldsymbol{D}}$.
    \begin{lemm}\label{lem:alpha_upperbd_underd}
        For all $\boldsymbol{v}\in \boldsymbol{Q}^\top(\mathcal{C}-\boldsymbol{x}^0)$ with $\|\boldsymbol{v}\|\le \min\{\appdiam(\boldsymbol{R}),1\}$, we have $\boldsymbol{v}=\boldsymbol{R}\boldsymbol{\alpha}_{\boldsymbol{v}}$ for some $\boldsymbol{\alpha}_{\boldsymbol{v}}\in\mathbb{R}^p$ with
        \begin{equation*}
            \left\|\boldsymbol{\alpha}_{\boldsymbol{v}}\right\| \le \min\{\appdiam(\boldsymbol{R}),1\}\left\|\boldsymbol{R}^{-1}\right\|_{\boldsymbol{x}^0,\mathcal{C},\boldsymbol{D}}.
        \end{equation*}
    \end{lemm}
    \begin{proof}
        Since $\boldsymbol{D}$ is full-column rank, $\boldsymbol{R}$ is invertible.  Thus,
        \begin{equation*}
            \left\|\boldsymbol{\alpha}_{\boldsymbol{v}}\right\| = \left\|\boldsymbol{R}^{-1}\boldsymbol{v}\right\| \le \max\limits_{\substack{\boldsymbol{w}\in \boldsymbol{Q}^\top(\mathcal{C}-\boldsymbol{x}^0)\\ \|\boldsymbol{w}\|\le \min\{\appdiam(\boldsymbol{R}),1\}}}\left\|\boldsymbol{R}^{-1}\boldsymbol{w}\right\| = \min\{\appdiam(\boldsymbol{R}),1\}\left\|\boldsymbol{R}^{-1}\right\|_{\boldsymbol{x}^0,\mathcal{C},\boldsymbol{D}}.
        \end{equation*}

    $\hfill\qed$
    \end{proof}

    Now we prove that $\widehat{m}(\widehat{\boldsymbol{s}})$ belongs to a class of $(\mathcal{C},\boldsymbol{Q})$-fully linear models of $f$ at $\boldsymbol{x}^0$ parametrized by $\appdiam(\boldsymbol{R})$.  We note that the constants $\kappa_{ef}$ and $\kappa_{eg}$ in Theorem \ref{thm:mhatfullylinear} depend on $\boldsymbol{x}^0$ and $\boldsymbol{D}$.  We will see later in Subsection \ref{subsec:modelacc} that, for the models constructed and used in $\myalg$, these dependencies can be removed.
    \begin{theo}\label{thm:mhatfullylinear}
        Suppose that $f\in C^{1+}$ with Lipschitz constant $L_{\nabla f}$.  Then, for all $\widehat{\boldsymbol{s}}\in \boldsymbol{Q}^\top(\mathcal{C}-\boldsymbol{x}^0)$ with $\|\widehat{\boldsymbol{s}}\|\le\appdiam(\boldsymbol{R})$, we have
        \begin{equation}\label{ineq:mhatfullylinear_graderr}
            \max\limits_{\substack{\boldsymbol{d}\in \boldsymbol{Q}^\top(\mathcal{C}-\boldsymbol{x}^0)\\ \|\boldsymbol{d}\|\le 1}}\left|\left(\boldsymbol{Q}^\top\nabla f(\boldsymbol{x}^0+\boldsymbol{Q}\widehat{\boldsymbol{s}})-\nabla\widehat{m}(\widehat{\boldsymbol{s}})\right)^\top \boldsymbol{d}\right|\le\left(\frac{1}{2}L_{\nabla f}\left(2+\sqrt{p}\left\|\widehat{\boldsymbol{R}}^{-1}\right\|_{\boldsymbol{x}^0,\mathcal{C},\boldsymbol{D}}\right)\right)\appdiam(\boldsymbol{R})
        \end{equation}
        and
        \begin{equation}\label{ineq:mhatfullylinear_fncerr}
            \left|f(\boldsymbol{x}^0+\boldsymbol{Q}\widehat{\boldsymbol{s}})-\widehat{m}(\widehat{\boldsymbol{s}})\right|\le\left(\frac{1}{2}L_{\nabla f}\left(1+\sqrt{p}\left\|\widehat{\boldsymbol{R}}^{-1}\right\|_{\boldsymbol{x}^0,\mathcal{C},\boldsymbol{D}}\right)\right)\appdiam(\boldsymbol{R})^2,
        \end{equation}
        where $\widehat{\boldsymbol{R}}=\boldsymbol{R}\slash\appdiam(\boldsymbol{R})$.
        
        In particular, $\widehat{m}(\widehat{\boldsymbol{s}})$ belongs to a class of $(\mathcal{C},\boldsymbol{Q})$-fully linear models of $f$ at $\boldsymbol{x}^0$ parametrized by $\appdiam(\boldsymbol{R})$, with constants
        \begin{equation*}
            \kappa_{ef} = \frac{1}{2}L_{\nabla f}\left(1+\sqrt{p}\left\|\widehat{\boldsymbol{R}}^{-1}\right\|_{\boldsymbol{x}^0,\mathcal{C},\boldsymbol{D}}\right)~~~\text{and}~~~\kappa_{eg} = \frac{1}{2}L_{\nabla f}\left(2+\sqrt{p}\left\|\widehat{\boldsymbol{R}}^{-1}\right\|_{\boldsymbol{x}^0,\mathcal{C},\boldsymbol{D}}\right).
        \end{equation*}
    \end{theo}
    \begin{proof}
        Following the proof of \cite[Theorem 9.5]{audet2017derivative}, using Equation \eqref{eq:GSG(R)multiplyR} and the fact that $\nabla \widehat{m}(\widehat{\boldsymbol{s}})=\nabla_S\widehat{f}(\boldsymbol{0}_p;\boldsymbol{R})$ for all $\widehat{\boldsymbol{s}}$, we have, for all $i=1,\ldots,p$,
        \begin{align*}
            \left|\boldsymbol{r}_i^\top\Big(\boldsymbol{Q}^\top\nabla f(\boldsymbol{x}^0)-\nabla\widehat{m}(\widehat{\boldsymbol{s}})\Big)\right| &= \left|\boldsymbol{r}_i^\top\nabla\widehat{f}(\boldsymbol{0}_p)-\boldsymbol{r}_i^\top\nabla_S\widehat{f}(\boldsymbol{0}_p;\boldsymbol{R})\right|\\
            &= \left|\boldsymbol{r}_i^\top\nabla\widehat{f}(\boldsymbol{0}_p)-e_i^\top\delta_{\widehat{f}}(\boldsymbol{0}_p;\boldsymbol{R})\right|\\
            &= \left|\boldsymbol{r}_i^\top\nabla\widehat{f}(\boldsymbol{0}_p)-\widehat{f}(\boldsymbol{r}_i)+\widehat{f}(\boldsymbol{0}_p)\right|\\
            &\le \frac{1}{2}L_{\nabla f}\appdiam(\boldsymbol{R})^2,
        \end{align*}
        where the last inequality is from, e.g., \cite[Lemma 9.4]{audet2017derivative}.  Therefore,
        \begin{equation*}
            \left\|\Big(\boldsymbol{Q}^\top\nabla f(\boldsymbol{x}^0)-\nabla\widehat{m}(\widehat{\boldsymbol{s}})\Big)^\top \boldsymbol{R}\right\| = \sqrt{\sum_{i=1}^p\left|\boldsymbol{r}_i^\top\Big(\boldsymbol{Q}^\top\nabla f(\boldsymbol{x}^0)-\nabla\widehat{m}(\widehat{\boldsymbol{s}})\Big)\right|^2} \le \frac{1}{2}\sqrt{p}L_{\nabla f}\appdiam(\boldsymbol{R})^2.
        \end{equation*}
        
        Let $\boldsymbol{d}\in \boldsymbol{Q}^\top(\mathcal{C}-\boldsymbol{x}^0)$ with $\|\boldsymbol{d}\|\le 1$.  Applying Lemma \ref{lem:alpha_upperbd_underd} with $\boldsymbol{v}=\min\{\appdiam(\boldsymbol{R}),1\}\boldsymbol{d}$, we get
        \begin{align*}
            \left|\Big(\boldsymbol{Q}^\top\nabla f(\boldsymbol{x}^0)-\nabla\widehat{m}(\widehat{\boldsymbol{s}})\Big)^\top \boldsymbol{d}\right| &= \left|\Big(\boldsymbol{Q}^\top\nabla f(\boldsymbol{x}^0)-\nabla\widehat{m}(\widehat{\boldsymbol{s}})\Big)^\top \frac{1}{\min\{\appdiam(\boldsymbol{R}),1\}}\boldsymbol{R}\boldsymbol{\alpha}_{\boldsymbol{v}}\right|\\
            &\le \frac{1}{\min\{\appdiam(\boldsymbol{R}),1\}}\left\|\Big(\boldsymbol{Q}^\top\nabla f(\boldsymbol{x}^0)-\nabla\widehat{m}(\widehat{\boldsymbol{s}})\Big)^\top \boldsymbol{R}\right\|\left\|\boldsymbol{\alpha}_{\boldsymbol{v}}\right\|\\
            &\le \frac{1}{2}\sqrt{p}L_{\nabla f}\appdiam(\boldsymbol{R})^2\left\|\boldsymbol{R}^{-1}\right\|_{\boldsymbol{x}^0,\mathcal{C},\boldsymbol{D}}\\
            &= \frac{1}{2}\sqrt{p}L_{\nabla f}\left\|\widehat{\boldsymbol{R}}^{-1}\right\|_{\boldsymbol{x}^0,\mathcal{C},\boldsymbol{D}}\appdiam(\boldsymbol{R}).
        \end{align*}
        
        From here, for all $\widehat{\boldsymbol{s}}\in \boldsymbol{Q}^\top(\mathcal{C}-\boldsymbol{x}^0)$ with $\|\widehat{\boldsymbol{s}}\|\le\appdiam(\boldsymbol{R})$, we obtain
        \begin{align*}
            &~~\left|\Big(\boldsymbol{Q}^\top\nabla f(\boldsymbol{x}^0+\boldsymbol{Q}\widehat{\boldsymbol{s}})-\nabla\widehat{m}(\widehat{\boldsymbol{s}})\Big)^\top \boldsymbol{d}\right|\\
            &= \left|\Big(\boldsymbol{Q}^\top\nabla f(\boldsymbol{x}^0)-\nabla\widehat{m}(\widehat{\boldsymbol{s}})+\boldsymbol{Q}^\top\nabla f(\boldsymbol{x}^0+\boldsymbol{Q}\widehat{\boldsymbol{s}})-\boldsymbol{Q}^\top\nabla f(\boldsymbol{x}^0)\Big)^\top \boldsymbol{d}\right|\\
            &\le \left|\Big(\boldsymbol{Q}^\top\nabla f(\boldsymbol{x}^0)-\nabla\widehat{m}(\widehat{\boldsymbol{s}})\Big)^\top \boldsymbol{d}\right| + \left|\Big(\boldsymbol{Q}^\top\nabla f(\boldsymbol{x}^0+\boldsymbol{Q}\widehat{\boldsymbol{s}})-\boldsymbol{Q}^\top\nabla f(\boldsymbol{x}^0)\Big)^\top \boldsymbol{d}\right|\\
            &\le \frac{1}{2}\sqrt{p}L_{\nabla f}\left\|\widehat{\boldsymbol{R}}^{-1}\right\|_{\boldsymbol{x}^0,\mathcal{C},\boldsymbol{D}}\appdiam(\boldsymbol{R}) + L_{\nabla f}\appdiam(\boldsymbol{R})\\
            &= \left(\frac{1}{2}L_{\nabla f}\left(2+\sqrt{p}\left\|\widehat{\boldsymbol{R}}^{-1}\right\|_{\boldsymbol{x}^0,\mathcal{C},\boldsymbol{D}}\right)\right)\appdiam(\boldsymbol{R}),        
        \end{align*}
        which proves Inequality \eqref{ineq:mhatfullylinear_graderr}.

        For Inequality \eqref{ineq:mhatfullylinear_fncerr}, we use \cite[Lemma 9.4]{audet2017derivative} again to get
        \begin{align*}
            \left|f(\boldsymbol{x}^0+\boldsymbol{Q}\widehat{\boldsymbol{s}})-\widehat{m}(\widehat{\boldsymbol{s}})\right| &= \left|\widehat{f}(\widehat{\boldsymbol{s}})-\widehat{m}(\widehat{\boldsymbol{s}})\right|\\
            &= \left|\widehat{f}(\widehat{\boldsymbol{s}})-\widehat{f}(\boldsymbol{0}_p)-\nabla\widehat{f}(\boldsymbol{0}_p)^\top\widehat{\boldsymbol{s}}+\nabla\widehat{f}(\boldsymbol{0}_p)^\top\widehat{\boldsymbol{s}}-\nabla\widehat{m}(\widehat{\boldsymbol{s}})^\top\widehat{\boldsymbol{s}}\right|\\
            &\le \frac{1}{2}L_{\nabla f}\left\|\widehat{\boldsymbol{s}}\right\|^2 + \left|\widehat{\boldsymbol{s}}^\top\Big(\nabla\widehat{f}(\boldsymbol{0}_p)-\nabla\widehat{m}(\widehat{\boldsymbol{s}})\Big)\right|\\
            &\le \frac{1}{2}L_{\nabla f}\appdiam(\boldsymbol{R})^2 + \left|\widehat{\boldsymbol{s}}^\top\Big(\nabla\widehat{f}(\boldsymbol{0}_p)-\nabla\widehat{m}(\widehat{\boldsymbol{s}})\Big)\right|\\
            &= \frac{1}{2}L_{\nabla f}\appdiam(\boldsymbol{R})^2 + \left|\Big(\boldsymbol{Q}^\top\nabla f(\boldsymbol{x}^0)-\nabla\widehat{m}(\widehat{\boldsymbol{s}})\Big)^\top\widehat{\boldsymbol{s}}\right|.
        \end{align*}
        Since $\boldsymbol{Q}^\top(\mathcal{C}-\boldsymbol{x}^0)$ is convex and contains $\boldsymbol{0}_p$, and $\|\widehat{\boldsymbol{s}}\|\le\appdiam(\boldsymbol{R})$, we have $\boldsymbol{v}=\widehat{\boldsymbol{s}}\slash\max\{\|\widehat{\boldsymbol{s}}\|,1\}\in \boldsymbol{Q}^\top(\mathcal{C}-\boldsymbol{x}^0)$ and $\|\boldsymbol{v}\|\le\min\{\appdiam(\boldsymbol{R}), 1\}$. Applying Lemma \ref{lem:alpha_upperbd_underd} with $\boldsymbol{v}=\widehat{\boldsymbol{s}}\slash\max\{\|\widehat{\boldsymbol{s}}\|,1\}$, we obtain
        \begin{align}\label{ineq:mhatfncerror}
        \begin{split}
            &~~\left|f(\boldsymbol{x}^0+\boldsymbol{Q}\widehat{\boldsymbol{s}})-\widehat{m}(\widehat{\boldsymbol{s}})\right|\\
            &\le \frac{1}{2}L_{\nabla f}\appdiam(\boldsymbol{R})^2 + \max\{\|\widehat{\boldsymbol{s}}\|,1\}\left|\Big(\boldsymbol{Q}^\top\nabla f(\boldsymbol{x}^0)-\nabla\widehat{m}(\widehat{\boldsymbol{s}})\Big)^\top \boldsymbol{R}\boldsymbol{\alpha}_{\boldsymbol{v}}\right|\\
            &\le \frac{1}{2}L_{\nabla f}\appdiam(\boldsymbol{R})^2 + \max\{\|\widehat{\boldsymbol{s}}\|,1\}\frac{1}{2}\sqrt{p}L_{\nabla f}\appdiam(\boldsymbol{R})^2\min\{\appdiam(\boldsymbol{R}),1\}\left\|\boldsymbol{R}^{-1}\right\|_{\boldsymbol{x}^0,\mathcal{C},\boldsymbol{D}}\\
            &= \frac{1}{2}L_{\nabla f}\appdiam(\boldsymbol{R})^2 + \frac{1}{2}\sqrt{p}L_{\nabla f}\left\|\widehat{\boldsymbol{R}}^{-1}\right\|_{\boldsymbol{x}^0,\mathcal{C},\boldsymbol{D}}\appdiam(\boldsymbol{R})\min\{\appdiam(\boldsymbol{R}),1\}\max\{\|\widehat{\boldsymbol{s}}\|,1\}.
        \end{split}
        \end{align}
        Noticing that 
        \begin{equation*}
            \min\{\appdiam(\boldsymbol{R}),1\}\max\{\|\widehat{\boldsymbol{s}}\|,1\}=
            \begin{cases}
                \|\widehat{\boldsymbol{s}}\|\le\appdiam(\boldsymbol{R}),&~\text{if}~\|\widehat{\boldsymbol{s}}\|\ge 1,\\
                \min\{\appdiam(\boldsymbol{R}),1\}\le\appdiam(\boldsymbol{R}),&~\text{if}~\|\widehat{\boldsymbol{s}}\|\le 1,
            \end{cases}
        \end{equation*}
        Inequality \eqref{ineq:mhatfncerror} gives
        \begin{align*}
            \left|f(\boldsymbol{x}^0+\boldsymbol{Q}\widehat{\boldsymbol{s}})-\widehat{m}(\widehat{\boldsymbol{s}})\right| &\le \frac{1}{2}L_{\nabla f}\appdiam(\boldsymbol{R})^2 + \frac{1}{2}\sqrt{p}L_{\nabla f}\left\|\widehat{\boldsymbol{R}}^{-1}\right\|_{\boldsymbol{x}^0,\mathcal{C},\boldsymbol{D}}\appdiam(\boldsymbol{R})^2\\ &=\left(\frac{1}{2}L_{\nabla f}\left(1+\sqrt{p}\left\|\widehat{\boldsymbol{R}}^{-1}\right\|_{\boldsymbol{x}^0,\mathcal{C},\boldsymbol{D}}\right)\right)\appdiam(\boldsymbol{R})^2.
        \end{align*}

    $\hfill\qed$
    \end{proof}

\section{Convex-constrained linear approximation random subspace trust-region algorithm}\label{sec:CLARSTA}
    In this section, we first introduce the framework of $\myalg$.  Then, we show that the models constructed and managed by $\myalg$ are $(\mathcal{C},\boldsymbol{Q})$-fully linear with constants $\kappa_{ef}$ and $\kappa_{eg}$ independent of $k$.  Last, we explain how the random subspaces are sampled and guaranteed to have a certain quality.
    
    $\myalg$ aims to find a first-order critical point.  In convex-constrained settings, one first-order criticality measure, introduced in \cite[Section 12.1.4]{conn2000trust}, is 
    \begin{equation*}
        \pi^f(\boldsymbol{x})=\left|\min\limits_{\substack{\boldsymbol{d}\in \mathcal{C}-\boldsymbol{x}\\ \|\boldsymbol{d}\|\le 1}}\nabla f(\boldsymbol{x})^\top \boldsymbol{d}\right|,~~~\text{for $\boldsymbol{x}\in \mathcal{C}$}.
    \end{equation*}
    Using the model $\widehat{m}$, we can define an approximate first-order criticality measure by
    \begin{equation*}
        \pi^m(\boldsymbol{x})=\left|\min\limits_{\substack{\boldsymbol{d}\in \boldsymbol{Q}^\top(\mathcal{C}-\boldsymbol{x})\\ \|\boldsymbol{d}\|\le 1}}\nabla \widehat{m}(\boldsymbol{0}_p)^\top \boldsymbol{d}\right|,~~~\text{for $\boldsymbol{x}\in \mathcal{C}$}.
    \end{equation*}

\subsection{Complete algorithm}
    In this subsection, we give the general framework of $\myalg$.  Similar to \cite{chen2024qfully}, it consists of the main algorithm (Algorithm~\ref{alg:mainAlg}) and three subroutines (Algorithms~\ref{alg:mainrm} to \ref{alg:dirngen}).  

    The sample directions in $\boldsymbol{D}_k$ managed by the general version of $\myalg$ consist of directions from some previously evaluated sample points, $\boldsymbol{D}_k^U$, and randomly generated new directions, $\boldsymbol{D}_k^R$.  At iterations outside $\mathcal{R}$, previously evaluated directions may be retained.  Here, $\mathcal{R}$ is a deterministic set of refresh iterations such that every block of $T$ consecutive iterations contains at least one index in $\mathcal{R}$.  When $\boldsymbol{D}_k^U$ is updated, its columns are selected by Algorithms~\ref{alg:mainrm} and \ref{alg:ptrm} in order to save the number of function evaluations while maintaining a reasonable geometry measure.  At each iteration in $\mathcal{R}$, all stored directions are discarded, and all $p$ directions are generated by Algorithm~\ref{alg:dirngen}.  Algorithm~\ref{alg:mainrm} includes an optional final step that allows the user to control the number of directions to be removed.  Algorithm~\ref{alg:ptrm} is inspired by \cite[Algorithm~4]{cartis2023scalable}, which itself is based on~\cite{powell2009bobyqa}.  In particular, the definition of $\theta_i$ follows the same philosophy: it penalizes points that are excessively distant while ensuring that the remaining matrix has a large minimum singular value.  Algorithms~\ref{alg:mainrm} and \ref{alg:ptrm} are used to maintain the radius and geometry conditions required for model accuracy; they are not intended to ensure that the subspace generated from reused directions is $\alpha$-well-aligned. The required probabilistic well-alignment guarantee is instead obtained at the prescribed full-refresh iterations, when all $p$ directions are generated afresh by Algorithm~\ref{alg:dirngen}.

    As a useful special case, we distinguish a minimal full-refresh version of $\myalg$ from the general sample-reuse version described above. In the minimal full-refresh version, all previously stored directions are discarded before constructing every model, and all $p$ directions are generated by Algorithm~\ref{alg:dirngen}.  Thus, $\boldsymbol{D}_k^U$ is empty and $\boldsymbol{D}_k=\boldsymbol{D}_k^R$ for every $k$, and the parameter $p_{\mathrm{rand}}$ is unnecessary; for compatibility with the input list of Algorithm~\ref{alg:mainAlg}, it may simply be set to $p$.  This version uses only Algorithms~\ref{alg:mainAlg} and \ref{alg:dirngen}.  It is obtained from the general version by taking $\mathcal{R}=\mathbb{N}$ and $T=1$. A high-level flowchart of this minimal full-refresh version of $\myalg$ can be found in Appendix \ref{app:flowchart}.
    
    In both the general sample-reuse version and the minimal full-refresh version, Algorithm~\ref{alg:dirngen} generates $\boldsymbol{D}_k^R$.  Its columns are mutually orthogonal and are also orthogonal to all columns of $\boldsymbol{D}_k^U$; in the minimal full-refresh version, $\boldsymbol{D}_k^U$ is empty. Indeed, when $\boldsymbol{Q}$ is specified as an orthogonal basis of $\mathrm{col}(\boldsymbol{D}_k^U)$, the operation $\widetilde{\boldsymbol{A}}=\boldsymbol{A}-\boldsymbol{Q}\boldsymbol{Q}^\top\boldsymbol{A}$ projects the columns of $\boldsymbol{A}$ onto $\mathrm{col}(\boldsymbol{D}_k^U)^\perp$, yielding the columns of $\widetilde{\boldsymbol{A}}$.  Thus, the columns of $\widetilde{\boldsymbol{Q}}$ are orthogonal to each other and are contained in $\mathrm{col}(\boldsymbol{D}_k^U)^\perp$.  This choice of new directions is optimal.  A detailed analysis of this sample set management procedure can be found in \cite[Subsection 2.3]{chen2024qfully}.  We also note that the requirement that $\boldsymbol{A}$ has full-column rank in Algorithm~\ref{alg:dirngen} is not restrictive, as i.i.d. $A_{ij}\sim\mathcal{N}(0,1)$ yields a full-column rank $\boldsymbol{A}$ with probability one (see, e.g., \cite[Proposition 7.1]{eaton2007multivariate}).

    We increase the flexibility of our algorithm by allowing $\boldsymbol{x}_{k+1}$ to be given by other optimization algorithms (e.g., direct-search methods or heuristics), as long as it satisfies Inequality \eqref{ineq:nextiterate}.  Note that the right-hand side of Inequality \eqref{ineq:nextiterate} is the minimum of all objective function values evaluated in the current iteration.  Moreover, instead of setting the trust-region radius increase parameter $\gamma_{\mathrm{inc}}$ to a fixed constant, we allow it to have different values in each iteration.  However, as shown in Section \ref{sec:converg}, we must set $\gamma_{\mathrm{inc}}^k = 1$ for all $k \ge \overline{k}$, for some $\overline{k} \ge 0$, in order to guarantee the convergence of $\myalg$.  We note that both $\{\gamma_{\mathrm{inc}}^k:k=0,1,\ldots\}$ and $\overline{k}$ are specified before running the algorithm and are never updated during its execution.
    
    It is worth mentioning that sample points in the form $\boldsymbol{x}_k+\boldsymbol{d}_i$ may not be feasible.  This means we require the constraints to be relaxable. Nevertheless, as we require $\boldsymbol{x}_0\in \mathcal{C}$ and only consider feasible points as candidates for the next iterate, all iterates $\boldsymbol{x}_k$ are feasible.  We note that the way we define~$\boldsymbol{s}_k$ (Equation \eqref{eq:sk}) is such that $\boldsymbol{x}_k+\boldsymbol{s}_k$ is the closest feasible point to the candidate solution $\boldsymbol{x}_k+\boldsymbol{Q}_k\widehat{\boldsymbol{s}}_k$ given by the subspace trust-region subproblem, so the right-hand side of Inequality~\eqref{ineq:nextiterate} is always well-defined.  For optimization problems where the constraints are not relaxable, we suggest trying $\boldsymbol{x}_k-\boldsymbol{d}_i$ for each $\boldsymbol{x}_k+\boldsymbol{d}_i\notin \mathcal{C}$.  If $\boldsymbol{x}_k-\boldsymbol{d}_i$ is also infeasible, then resample the $\boldsymbol{d}_i$.  While this procedure works well in practice, theoretical analyses are relegated to future work.

    \begin{algorithm}[htbp!]\caption{Main algorithm}\label{alg:mainAlg}
    \DontPrintSemicolon
        \KwIn{\begin{tabular}{lcl}
            $f$ & $f:\mathbb{R}^n\to\mathbb{R}$ & objective function\\
            $p$ & $1\le p\le n$ & subspace dimension\\
            $p_{\mathrm{rand}}$ & $1\le p_{\mathrm{rand}} \le p$ & minimum randomized subspace dimension\\
            $\boldsymbol{x}_0$ & $\boldsymbol{x}_0\in \mathcal{C}$ & feasible starting point\\
            $\Delta_{\min}$, $\Delta_{\max}$ &  $0 < \Delta_{\min} \le \Delta_{\max}$ & minimum and maximum trust-region radii\\
            $\Delta_0$ & $\Delta_{\min} \le \Delta_0 \le \Delta_{\max}$ & initial trust-region radius\\
            $\gamma_{\mathrm{dec}}$ & $0 < \gamma_{\mathrm{dec}} < 1$ & trust-region decrease scaling\\
            $\gamma_{\mathrm{inc}}^k$ & $\gamma_{\mathrm{inc}}^k \ge 1, k\in\mathbb{N}$ & per-iteration trust-region increase scaling\\
            $\eta_1$, $\eta_2$ & $0 < \eta_1 \le \eta_2 < 1$ & acceptance thresholds\\
            $\mu$ & $\mu>0$ & criticality constant\\
            $\epsilon_{\mathrm{rad}}$, $\epsilon_{\mathrm{geo}}$ & $\epsilon_{\mathrm{rad}}\ge 1$, $\epsilon_{\mathrm{geo}}>0$ & sample set radius and geometry tolerance\\
            $\mathcal{R}$ & $\mathcal{R}\subseteq\mathbb{N}$ & prescribed deterministic refresh schedule
        \end{tabular}}
        {\bf Initialize:} Generate random mutually orthogonal directions $\boldsymbol{d}_1,\ldots, \boldsymbol{d}_p\in\mathbb{R}^n$ such that each $\|\boldsymbol{d}_i\|\le\Delta_0$ by Algorithm~\ref{alg:dirngen}; define initial direction matrix $\boldsymbol{D}_0=[\boldsymbol{d}_1\cdots \boldsymbol{d}_p]$.
        
        \For{$k=0,1,\ldots$}
        { 
            Perform the $\boldsymbol{Q}\boldsymbol{R}$-factorization $\boldsymbol{D}_k=\boldsymbol{Q}_k\boldsymbol{R}_k$ and construct model $\widehat{m}_k$ at $\boldsymbol{x}_k$.
            
            \uIf{$\Delta_k\le\mu\pi^m(\boldsymbol{x}_k)$ \tcp*[f]{Model accuracy test}}
            {
                \tcp{Trust-region subproblem and update}
                Solve (approximately)
                \begin{equation*}
                    \widehat{\boldsymbol{s}}_k \in \argmin\left\{\widehat{m}_k(\widehat{\boldsymbol{s}}):\widehat{\boldsymbol{s}}\in \boldsymbol{Q}_k^\top\left(\mathcal{C}-\boldsymbol{x}_k\right), \left\|\widehat{\boldsymbol{s}}\right\|\le\Delta_k\right\}
                \end{equation*}
                and let
                \begin{equation}\label{eq:sk}
                    \boldsymbol{s}_k=\mathrm{proj}_{\mathcal{C}}(\boldsymbol{x}_k+\boldsymbol{Q}_k\widehat{\boldsymbol{s}}_k)-\boldsymbol{x}_k.
                \end{equation}
                        
                Calculate
                \begin{equation*}
                    \rho_k=\frac{f(\boldsymbol{x}_k)-f(\boldsymbol{x}_k+\boldsymbol{s}_k)}{\widehat{m}_k(\boldsymbol{0}_p)-\widehat{m}_k(\widehat{\boldsymbol{s}}_k)}
                \end{equation*}
                and update the trust-region radius
                \begin{equation*}
                  \Delta_{k+1}=
                    \begin{cases}
                      \gamma_{\mathrm{dec}}\Delta_k, & \text{if }\rho_k<\eta_1,\\
                      \min\left(\gamma_{\mathrm{inc}}^k\Delta_k, \Delta_{\max}\right), &\text{if } \rho_k>\eta_2,\\
                      \Delta_k, & \text{otherwise}.
                    \end{cases}       
                \end{equation*}
                
                Let $\boldsymbol{x}_{k+1}$ be any point such that $\boldsymbol{x}_{k+1}\in \mathcal{C}$ and
                \begin{equation}\label{ineq:nextiterate}
                    f(\boldsymbol{x}_{k+1})\le\min\left\{\{f(\boldsymbol{x}_k+\boldsymbol{s}_k)\}\cup\{f(\boldsymbol{x}_k+\boldsymbol{d}_i):\boldsymbol{x}_k+\boldsymbol{d}_i\in \mathcal{C}, \boldsymbol{d}_i\in [\boldsymbol{0}_n~\boldsymbol{D}_k]\}\right\}.   
                \end{equation}

                If $k+1\in\mathcal{R}$, set $\boldsymbol{D}^U_{k+1}$ to be empty; otherwise, create $\boldsymbol{D}^U_{k+1}$ according to Algorithms~\ref{alg:mainrm} and \ref{alg:ptrm}.  Let $p_1$ be the number of columns of $\boldsymbol{D}^U_{k+1}$.
                
                Use Algorithm~\ref{alg:dirngen} with $q=p-p_1$ and $\Delta=\Delta_{k+1}$ to generate $\boldsymbol{D}_{k+1}^R$ orthogonal to $\boldsymbol{D}_{k+1}^U$, and set $\boldsymbol{D}_{k+1}=[\boldsymbol{D}_{k+1}^U~\boldsymbol{D}_{k+1}^R]$.
            }
            \Else
            {
                Set $\Delta_{k+1}=\gamma_{\mathrm{dec}}\Delta_k$ and $\boldsymbol{x}_{k+1}=\boldsymbol{x}_k$.  If $k+1\in\mathcal{R}$, set $\boldsymbol{D}_{k+1}^U$ to be empty and generate $\boldsymbol{D}_{k+1}=\boldsymbol{D}_{k+1}^R$ by Algorithm~\ref{alg:dirngen} with $q=p$ and $\Delta=\Delta_{k+1}$; otherwise, set $\boldsymbol{D}_{k+1}=\gamma_{\mathrm{dec}}\boldsymbol{D}_k$.
            }
            \lIf{$\Delta_{k+1}<\Delta_{\min}$,}{\textbf{stop} \tcp*[f]{Stopping test}}
        }
    \end{algorithm}

    \begin{algorithm}[htbp!]\caption{Sample-reuse $\boldsymbol{D}^U_{k+1}$ algorithm (inspired by \cite[Algorithm~2]{chen2024qfully})}\label{alg:mainrm}
    \DontPrintSemicolon
        \KwIn{$\boldsymbol{D}_k$; $\boldsymbol{x}_k$; $\boldsymbol{x}_{k+1}$; $\boldsymbol{s}_k$; $\Delta_{k+1}$; $p$; $p_{\mathrm{rand}}$; $\epsilon_{\mathrm{rad}}$; $\epsilon_{\mathrm{geo}}$.}
        Select $p$ linearly independent directions from $\{\{\boldsymbol{x}_k-\boldsymbol{x}_{k+1}\}\cup\{\boldsymbol{x}_k+\boldsymbol{s}_k-\boldsymbol{x}_{k+1}\}\cup\{\boldsymbol{x}_k+\boldsymbol{d}_i-\boldsymbol{x}_{k+1}:\boldsymbol{d}_i\in [\boldsymbol{0}_n~\boldsymbol{D}_k]\}\}$ as the columns of $\boldsymbol{D}^U_{k+1}$.
        
        Remove $p_{\mathrm{rand}}$ directions from $\boldsymbol{D}^U_{k+1}$ using Algorithm~\ref{alg:ptrm}.
        
        Remove all directions $\boldsymbol{d}\in \boldsymbol{D}^U_{k+1}$ with $\|\boldsymbol{d}\|>\epsilon_{\mathrm{rad}}\Delta_{k+1}$.
        
        \While{$\boldsymbol{D}^U_{k+1}\text{ is not empty and }\sigma_{\min}(\boldsymbol{D}^U_{k+1})<\epsilon_{\mathrm{geo}}$}
        {
            Remove 1 direction from $\boldsymbol{D}^U_{k+1}$ using Algorithm~\ref{alg:ptrm}.
        }

       {\em (optional)}  If $\boldsymbol{D}^U_{k+1}$ is not empty and the above steps terminate before the desired number of removals is reached, continue removing directions using Algorithm~\ref{alg:ptrm}. 
    \end{algorithm}
    
    \begin{algorithm}[htbp!]\caption{Direction removing algorithm (from \cite[Algorithm~3]{chen2024qfully})}\label{alg:ptrm}
    \DontPrintSemicolon
        \KwIn{$\boldsymbol{D}^U_{k+1}$; $\Delta_{k+1}$; number of directions to be removed $p_{\mathrm{rm}}$.}
        Set the number of directions removed {\tt removed} = 0.
        
        \While{${\tt removed}<p_{\mathrm{rm}}$}
        {
            Denote $\boldsymbol{D}^U_{k+1}=[\boldsymbol{d}^U_1\cdots \boldsymbol{d}^U_m]$.
            
            \For{$i=1,\ldots,m$}
            {
                Define matrix $\boldsymbol{M}_i=\left[\boldsymbol{d}^U_1 \cdots \boldsymbol{d}^U_{i-1} ~\boldsymbol{d}^U_{i+1} \cdots \boldsymbol{d}^U_m\right]$ and compute 
                \begin{equation*}
                    \theta_i=\sigma_{\min}(\boldsymbol{M}_i)\cdot\max\left(\frac{\left\|\boldsymbol{d}^U_i\right\|^4}{\Delta^4_{k+1}},1\right).
                \end{equation*}
            }
            Remove the direction with the largest $\theta_i$ from $\boldsymbol{D}^U_{k+1}$, set ${\tt removed}={\tt removed}+1$.
        }
    \end{algorithm}

    \begin{algorithm}[htbp!]\caption{Direction generating algorithm (inspired by \cite[Algorithm~4.1]{hough2022model})}\label{alg:dirngen}
    \DontPrintSemicolon
        \KwIn{Orthogonal basis for current subspace $\boldsymbol{Q}$ (optional); number of new directions $q$; length of new directions upper bound $\Delta>0$.}
        Generate $\boldsymbol{A}\in\mathbb{R}^{n\times q}$ with i.i.d. $A_{ij}\sim\mathcal{N}(0,1)$ such that $\boldsymbol{A}$ has full-column rank.
        
        If $\boldsymbol{Q}$ is specified, then calculate $\widetilde{\boldsymbol{A}}=\boldsymbol{A}-\boldsymbol{Q}\boldsymbol{Q}^\top \boldsymbol{A}$, otherwise set $\widetilde{\boldsymbol{A}}=\boldsymbol{A}$.
    
        Perform the $\boldsymbol{Q}\boldsymbol{R}$-factorization $\widetilde{\boldsymbol{A}}=\widetilde{\boldsymbol{Q}}\widetilde{\boldsymbol{R}}$ and denote $\widetilde{\boldsymbol{Q}}=[\widetilde{\boldsymbol{q}}_1\cdots \widetilde{\boldsymbol{q}}_q]$.
        
        Return $\Delta\widetilde{\boldsymbol{q}}_1,\ldots,\Delta\widetilde{\boldsymbol{q}}_q$.
    \end{algorithm}\FloatBarrier

\subsection{Model accuracy}\label{subsec:modelacc}
    In this subsection, we prove that each $\widehat{m}_k$ belongs to a class of $(\mathcal{C},\boldsymbol{Q}_k)$-fully linear models of $f$ at $\boldsymbol{x}_k$ parametrized by $\Delta_k$, and the constants $\kappa_{ef}$ and $\kappa_{eg}$ are independent of $k$.  This is essential for the convergence of $\myalg$.

    In order to apply Theorem \ref{thm:mhatfullylinear} to reach our goal, we need to find the relation between $\appdiam(\boldsymbol{R}_k)$ and $\Delta_k$, and find an upper bound for $\|\widehat{\boldsymbol{R}}_k^{-1}\|_{\boldsymbol{x}_k,\mathcal{C},\boldsymbol{D}_k}$ that is independent of $k$.  
    
    First, from Algorithms~\ref{alg:mainrm} and \ref{alg:dirngen} we have all columns in $\boldsymbol{D}_k^U$ are bounded by $\epsilon_{\mathrm{rad}}\Delta_k$ and all columns in $\boldsymbol{D}_k^R$ have length $\Delta_k$.  Moreover, $\boldsymbol{D}_k^R$ is never empty, i.e., each $\boldsymbol{D}_k$ has at least one column generated by Algorithm~\ref{alg:dirngen}.  These imply that
    \begin{equation*}
        \Delta_k \le \appdiam(\boldsymbol{R}_k)=\max\limits_{1\le i\le p}\|\boldsymbol{r}_i\|=\max\limits_{1\le i\le p}\|\boldsymbol{Q}\boldsymbol{r}_i\|=\max\limits_{1\le i\le p}\|\boldsymbol{d}_i\|\le\epsilon_{\mathrm{rad}}\Delta_k.
    \end{equation*}
    Notice that $\|\widehat{\boldsymbol{R}}_k^{-1}\|_{\boldsymbol{x}_k,\mathcal{C},\boldsymbol{D}_k}\le\|\widehat{\boldsymbol{R}}_k^{-1}\|$ and an upper bound of $\|\widehat{\boldsymbol{R}}_k^{-1}\|$ that is independent of $k$ is 
    \begin{equation*}
        M_{\widehat{\boldsymbol{R}}^{-1}}=\max\left\{\frac{1}{\epsilon_{\mathrm{geo}}}, \frac{1}{\Delta_{\min}}\right\}\epsilon_{\mathrm{rad}}\Delta_{\max},
    \end{equation*}
    which is proved in \cite[Theorem 6 and Inequality (8)]{chen2024qfully}.

    Therefore, we have the following theorem, which shows that each $\widehat{m}_k$ belongs to a class of $(\mathcal{C},\boldsymbol{Q}_k)$-fully linear models of $f$ at $\boldsymbol{x}_k$ parametrized by $\Delta_k$, and the constants $\kappa_{ef}$ and $\kappa_{eg}$ are independent of $k$.
    \begin{theo}
        Suppose that $f\in C^{1+}$ with Lipschitz constant $L_{\nabla f}$.  Then, each $\widehat{m}_k$ belongs to a class of $(\mathcal{C},\boldsymbol{Q}_k)$-fully linear models of $f$ at $\boldsymbol{x}_k$ parameterized by $\Delta_k\in(0,\Delta_{\max}]$ with the following constants
        \begin{equation*}
            \kappa_{ef} = \left(\frac{1}{2}L_{\nabla f}\left(1+\sqrt{p}M_{\widehat{\boldsymbol{R}}^{-1}}\right)\right)\epsilon^2_{\mathrm{rad}}~~~\text{and}~~~\kappa_{eg} = \left(\frac{1}{2}L_{\nabla f}\left(2+\sqrt{p}M_{\widehat{\boldsymbol{R}}^{-1}}\right)\right)\epsilon_{\mathrm{rad}}.
        \end{equation*}
    \end{theo}
    \begin{proof}
          Taking $\boldsymbol{x}^0=\boldsymbol{x}_k$, $\boldsymbol{D}=\boldsymbol{D}_k$, $\boldsymbol{Q}=\boldsymbol{Q}_k$, $\boldsymbol{R}=\boldsymbol{R}_k$, and $\widehat{m}=\widehat{m}_k$ in Theorem~\ref{thm:mhatfullylinear}, using the fact that $\appdiam(\boldsymbol{R}_k)\le\epsilon_{\mathrm{rad}}\Delta_k$, we have that, for all $\widehat{\boldsymbol{s}}\in \boldsymbol{Q}_k^\top(\mathcal{C}-\boldsymbol{x}_k)$ with $\|\widehat{\boldsymbol{s}}\|\le\appdiam(\boldsymbol{R}_k)$,
        \begin{align*}
            \left|f(\boldsymbol{x}_k+\boldsymbol{Q}_k\widehat{\boldsymbol{s}})-\widehat{m}_k(\widehat{\boldsymbol{s}})\right| &\le \left(\frac{1}{2}L_{\nabla f}\left(1+\sqrt{p}\left\|\widehat{\boldsymbol{R}}_k^{-1}\right\|_{\boldsymbol{x}_k,\mathcal{C},\boldsymbol{D}_k}\right)\right)\epsilon_{\mathrm{rad}}^2\Delta_k^2\\
            \max\limits_{\substack{\boldsymbol{d}\in \boldsymbol{Q}_k^\top(\mathcal{C}-\boldsymbol{x}_k)\\ \|\boldsymbol{d}\|\le 1}}\left|\left(\boldsymbol{Q}_k^\top\nabla f(\boldsymbol{x}_k+\boldsymbol{Q}_k\widehat{\boldsymbol{s}})-\nabla\widehat{m}_k(\widehat{\boldsymbol{s}})\right)^\top \boldsymbol{d}\right| &\le \left(\frac{1}{2}L_{\nabla f}\left(2+\sqrt{p}\left\|\widehat{\boldsymbol{R}}_k^{-1}\right\|_{\boldsymbol{x}_k,\mathcal{C},\boldsymbol{D}_k}\right)\right)\epsilon_{\mathrm{rad}}\Delta_k.
        \end{align*}
        The results follow from using the fact that $\Delta_k\le\appdiam(\boldsymbol{R}_k)$ and $\|\widehat{\boldsymbol{R}}_k^{-1}\|_{\boldsymbol{x}_k,\mathcal{C},\boldsymbol{D}_k}\le M_{\widehat{\boldsymbol{R}}^{-1}}$.

    $\hfill\qed$
    \end{proof}

\subsection{Subspace quality}
    Other than model accuracy, the convergence of random subspace model-based DFO methods normally also requires the subspaces to have a certain level of quality.  For unconstrained problems, subspaces with such quality can be constructed using $\alpha$-well-aligned matrices \cite{cartis2023scalable} with the sampling methods provided in \cite{cartis2022randomiseda,cartis2022randomisedb,dzahini2024stochastic,shao2021random}.  
    
    In this subsection, we first define a normalized convex-constrained version of $\alpha$-well-aligned matrices.  Then, we provide a method for sampling these matrices.

    The $\alpha$-well-aligned matrices \cite{cartis2023scalable} preserve the unconstrained first-order criticality measure (i.e., $\|\nabla f\|$) by a certain fraction with certain probability lower bounds.  For convex-constrained problems, we can define matrices with similar properties.  We use the normalization $n\slash z$ to compensate for the average shrinkage caused by a rank-$z$ orthogonal projection matrix.
    \begin{defini}\label{def:generalizedAlphaWellAligned}
        Let $f\in C^1$, $\boldsymbol{x}\in\mathcal{C}$, and $\alpha\in (0,1)$.  Let $\boldsymbol{A}\in\mathbb{R}^{n\times z}$ have full-column rank.  Suppose that $\boldsymbol{A}=\boldsymbol{Q}\boldsymbol{R}$ is the $\boldsymbol{Q}\boldsymbol{R}$-factorization of $\boldsymbol{A}$.  We say that $\boldsymbol{A}$ is $\alpha$-well-aligned for $f$ and $\mathcal{C}$ at $\boldsymbol{x}$ if 
        \begin{equation}\label{ineq:generalizedAlphaWellAligned}
            \left|\min\limits_{\substack{\boldsymbol{d}\in \mathcal{C}-\boldsymbol{x}\\ \|\boldsymbol{d}\|\le 1}}\nabla f(\boldsymbol{x})^\top \left(\frac{n}{z}\boldsymbol{Q}\boldsymbol{Q}^\top\right) \boldsymbol{d}\right| \ge \alpha\pi^f(\boldsymbol{x}).
        \end{equation}
    \end{defini}
    \begin{rema}
        If $\mathcal{C}=\mathbb{R}^n$, then Inequality \eqref{ineq:generalizedAlphaWellAligned} becomes
        \begin{equation*}
             \frac{n}{z}\left\|\boldsymbol{Q}\boldsymbol{Q}^\top\nabla f(\boldsymbol{x})\right\| = \frac{n}{z}\left\|\boldsymbol{Q}^\top\nabla f(\boldsymbol{x})\right\| \ge \alpha\left\|\nabla f(\boldsymbol{x})\right\|.
        \end{equation*}
        Equivalently, $\boldsymbol{A}$ is a standard $(\alpha z\slash n)$-well-aligned matrix in the unconstrained sense of \cite{cartis2023scalable}.
    \end{rema}
    \begin{rema}
        Similar to the definition of $\alpha$-well-aligned matrices in the unconstrained case \cite{cartis2023scalable}, the left-hand side of Inequality \eqref{ineq:generalizedAlphaWellAligned} is given by replacing $\nabla f(\boldsymbol{x})$ in the first-order criticality measure with a projected gradient.  The difference is the normalization $n\slash z$.  If $\boldsymbol{X}\in\mathbb{R}^{n\times n}$ is a uniformly distributed orthogonal projection matrix of rank $z$, then $\mathbb{E}[(n\slash z)\boldsymbol{X}]=\boldsymbol{I}_n$.  Therefore, Inequality~\eqref{ineq:generalizedAlphaWellAligned} can be interpreted as requiring the normalized projected gradient to preserve at least a certain fraction of the first-order criticality measure.
    \end{rema}
    
    The following lemma provides insight into how to sample $\alpha$-well-aligned matrices.  In particular, suppose that $\boldsymbol{A}$ is sampled so that $\boldsymbol{Q}\boldsymbol{Q}^\top$ satisfies the condition on $\boldsymbol{X}$ stated in the lemma.  Then, by taking $\boldsymbol{v}=\nabla f(\boldsymbol{x})$, the lemma guarantees that $\boldsymbol{A}$ is $\alpha$-well-aligned for $f$ and $\mathcal{C}$ at $\boldsymbol{x}$ with a probability lower bound independent of $n$.
    \begin{lemm}\label{lem:generalizedAlphaWellAligned_prob}
        Let $\boldsymbol{x}\in\mathcal{C}$, $\alpha\in(0,1)$, and integer $z\in(0,n]$.  Let $P_{z,n}$ denote the set of all $n\times n$ orthogonal projection matrices of rank $z$.  Suppose that $\boldsymbol{X}\in\mathbb{R}^{n\times n}$ is uniformly distributed on~$P_{z,n}$.  Then, for any deterministic vector $\boldsymbol{v}\in\mathbb{R}^n$,
        \begin{equation*}
            \mathbb{P}\left[\left|\min\limits_{\substack{\boldsymbol{d}\in \mathcal{C}-\boldsymbol{x}\\ \|\boldsymbol{d}\|\le 1}}\boldsymbol{v}^\top \left(\frac{n}{z}\boldsymbol{X}\right) \boldsymbol{d}\right| \ge \alpha \left|\min\limits_{\substack{\boldsymbol{d}\in \mathcal{C}-\boldsymbol{x}\\ \|\boldsymbol{d}\|\le 1}}\boldsymbol{v}^\top \boldsymbol{d}\right|\right]\ge r_{\alpha,z},~~~\text{where}~~~r_{\alpha,z}=\frac{z(1-\alpha)^2}{2(z+2)}>0.
        \end{equation*}
        If
        \begin{equation*}
            \left|\min\limits_{\substack{\boldsymbol{d}\in \mathcal{C}-\boldsymbol{x}\\ \|\boldsymbol{d}\|\le 1}}\boldsymbol{v}^\top \boldsymbol{d}\right|=0,
        \end{equation*}
        then the event in the probability holds with probability one.
    \end{lemm}
    \begin{proof}
        Let
        \begin{equation*}
            \boldsymbol{d}^*\in\argmin\left\{\boldsymbol{v}^\top \boldsymbol{d}:\boldsymbol{d}\in \mathcal{C}-\boldsymbol{x}, \|\boldsymbol{d}\|\le 1\right\}.
        \end{equation*}
        Notice that $\boldsymbol{0}_n\in\mathcal{C}-\boldsymbol{x}$ and $\|\boldsymbol{0}_n\|\le 1$, so we must have $\boldsymbol{v}^\top\boldsymbol{d}^*\le 0$ and 
        \begin{equation*}
            \left|\min\limits_{\substack{\boldsymbol{d}\in \mathcal{C}-\boldsymbol{x}\\ \|\boldsymbol{d}\|\le 1}}\boldsymbol{v}^\top \boldsymbol{d}\right| = -\min\limits_{\substack{\boldsymbol{d}\in \mathcal{C}-\boldsymbol{x}\\ \|\boldsymbol{d}\|\le 1}}\boldsymbol{v}^\top \boldsymbol{d} = -\boldsymbol{v}^\top \boldsymbol{d}^*.
        \end{equation*}
        When $-\boldsymbol{v}^\top \boldsymbol{d}^*=0$, the result is obvious.  Suppose that $-\boldsymbol{v}^\top\boldsymbol{d}^*>0$.  Define
        \begin{equation*}
            \boldsymbol{u}=-\frac{\boldsymbol{v}}{\|\boldsymbol{v}\|},~~~\boldsymbol{w}=\frac{\boldsymbol{d}^*}{\|\boldsymbol{d}^*\|},~~~\text{and}~~~\theta=\boldsymbol{u}^\top\boldsymbol{w}=\frac{-\boldsymbol{v}^\top\boldsymbol{d}^*}{\|\boldsymbol{v}\|\|\boldsymbol{d}^*\|}\in(0,1].
        \end{equation*}
        Define $Y=\boldsymbol{u}^\top\boldsymbol{X}\boldsymbol{u}$.  If $\theta<1$, let $\boldsymbol{s}$ be the unit vector orthogonal to $\boldsymbol{u}$ such that
        \begin{equation*}
            \boldsymbol{w}=\theta\boldsymbol{u}+\sqrt{1-\theta^2}\boldsymbol{s},
        \end{equation*}
        and define $Z=\boldsymbol{u}^\top\boldsymbol{X}\boldsymbol{s}$.  In this case, on the event $\{Y\ge\alpha z\slash n, Z\ge0\}$, we have
        \begin{align*}
            -\boldsymbol{v}^\top\left(\frac{n}{z}\boldsymbol{X}\right)\boldsymbol{d}^* &= \frac{n}{z}\|\boldsymbol{v}\|\|\boldsymbol{d}^*\|\boldsymbol{u}^\top\boldsymbol{X}\boldsymbol{w}\\
            &=\frac{n}{z}\|\boldsymbol{v}\|\|\boldsymbol{d}^*\|\left(\theta Y+\sqrt{1-\theta^2}Z\right)\\
            &\ge\alpha\theta\|\boldsymbol{v}\|\|\boldsymbol{d}^*\|\\
            &=\alpha\left(-\boldsymbol{v}^\top\boldsymbol{d}^*\right).
        \end{align*}
        Since $\boldsymbol{d}^*$ is feasible in the minimization, this implies the event in the statement of the lemma.  If $\theta=1$, the same implication holds on the event $\{Y\ge\alpha z\slash n\}$ without introducing $Z$.

        It remains to bound the probability of this event.  When $\theta<1$, let $\boldsymbol{H}=\boldsymbol{I}_n-2\boldsymbol{s}\boldsymbol{s}^\top$.  Since $\|\boldsymbol{s}\|=1$, we have that $\boldsymbol{H}$ is orthogonal.  By the orthogonal invariance of the uniform distribution on $P_{z,n}$ \cite[Theorem 2.2.2]{chikuse2012statistics}, $\boldsymbol{H}\boldsymbol{X}\boldsymbol{H}=\boldsymbol{H}\boldsymbol{X}\boldsymbol{H}^\top$ has the same distribution as $\boldsymbol{X}$.  Moreover, this transformation leaves $Y$ unchanged and changes $Z$ to $-Z$.  Hence,
        \begin{equation}\label{eq:symmetry_YZ}
            \mathbb{P}\left[Y\ge\alpha\frac{z}{n},Z\ge0\right]\ge\frac{1}{2}\mathbb{P}\left[Y\ge\alpha\frac{z}{n}\right].
        \end{equation}
        By \cite[Theorem 2.2.2]{chikuse2012statistics}, if $\boldsymbol{Q}$ is uniformly distributed on the set of $n\times n$ orthonormal matrices $O(n)$, then $\boldsymbol{X}$ has the same distribution as $$\boldsymbol{Q}\begin{bmatrix}
            \boldsymbol{I}_z & \boldsymbol{0}_{z\times(n-z)}\\
            \boldsymbol{0}_{(n-z)\times z} & \boldsymbol{0}_{(n-z)\times(n-z)}
        \end{bmatrix}\boldsymbol{Q}^{\top}.$$  For a fixed unit vector $\boldsymbol{u}$, let
        $\boldsymbol{\xi}=\boldsymbol{Q}^{\top}\boldsymbol{u}$. By the
        orthogonal invariance of the uniform distribution on $O(n)$
        \cite[Theorem~2.2.1]{chikuse2012statistics}, $\boldsymbol{\xi}$ is uniformly
        distributed on the unit sphere. Therefore, $Y$ has the same distribution as $\sum_{i=1}^{z}\xi_i^2$. Using the fact that $\mathbb{E}[\xi_i^2]=1\slash n$, $\mathbb{E}[\xi_i^4]=3\slash(n(n+2))$, and $\mathbb{E}[\xi_i^2\xi_j^2]=1\slash(n(n+2))$ for $i\neq j$ \cite[Proposition 2.5]{meckes2019random}, we obtain
        \begin{equation*}
            \mathbb{E}[Y]=\frac{z}{n}~~~\text{and}~~~\mathbb{E}[Y^2]=\frac{z(z+2)}{n(n+2)}.
        \end{equation*}
        By the Paley-Zygmund inequality (e.g., \cite[Exercise 1.16]{Vershynin_2026}),
        \begin{align*}
            \mathbb{P}\left[Y\ge\alpha\frac{z}{n}\right] \ge(1-\alpha)^2\frac{\mathbb{E}[Y]^2}{\mathbb{E}[Y^2]} =(1-\alpha)^2\frac{z(n+2)}{n(z+2)} \ge(1-\alpha)^2\frac{z}{z+2}.
        \end{align*}
        Combining this bound with Inequality \eqref{eq:symmetry_YZ} proves the result when $\theta<1$.  When $\theta=1$, the factor $1\slash2$ is not needed, so the stated lower bound also holds.

    $\hfill\qed$
    \end{proof}
    \begin{rema}
        The uniform distribution on $P_{z,n}$ is the unique probability measure invariant under orthogonal conjugation. It is obtained by pushing forward the normalized Haar measure on $O(n)$.  We refer the reader to \cite{chikuse2012statistics} for details. 
    \end{rema}

    Using Lemma \ref{lem:generalizedAlphaWellAligned_prob}, we can prove that, whenever $\boldsymbol{D}_k^U$ is empty, the matrix $\boldsymbol{D}_k$ constructed in $\myalg$ is $\alpha$-well-aligned for $f$ and $\mathcal{C}$ at $\boldsymbol{x}_k$ for any fixed $\alpha\in(0,1)$, with a probability lower bound independent of $n$.
    \begin{theo}\label{thm:gen_alphawellaligned}
        Suppose that $f\in C^1$.  Let $\alpha\in(0,1)$ and integer $p\in(0,n]$.  Suppose that $\boldsymbol{D}_k^U$ is empty and~$\boldsymbol{D}_k^R$ is constructed by Algorithm~\ref{alg:dirngen} with $\Delta=\Delta_k$.  Conditional on $\boldsymbol{x}_k$ and all information available immediately before the random directions forming $\boldsymbol{D}_k$ are generated, $\boldsymbol{D}_k=\boldsymbol{D}_k^R\in\mathbb{R}^{n\times p}$ is $\alpha$-well-aligned for $f$ and~$\mathcal{C}$ at $\boldsymbol{x}_k$ with probability at least
        \begin{equation*}
            r_{\alpha,p}=\frac{p(1-\alpha)^2}{2(p+2)}>0.
        \end{equation*}
        In particular, both $\alpha$ and $p$ may be selected independently of $n$.
    \end{theo}
    \begin{proof}
        Let $\boldsymbol{A}$ be the matrix in Algorithm~\ref{alg:dirngen} with $A_{ij}\sim\mathcal{N}(0,1)$.   Since $p\le n$, $\boldsymbol{A}$ has full-column rank with probability one.  Moreover, since $\boldsymbol{D}_k^U$ is empty, we have $\widetilde{\boldsymbol{A}}=\boldsymbol{A}$ in Algorithm~\ref{alg:dirngen}, $p_1=0$, and $q=p$ in Algorithm~\ref{alg:mainAlg}.
        
        Let $\boldsymbol{D}_k=\boldsymbol{Q}_k\boldsymbol{R}_k$ be the $\boldsymbol{Q}\boldsymbol{R}$-factorization of $\boldsymbol{D}_k$.  According to \cite{muirhead2009aspects}, $\boldsymbol{Q}_k$ is uniformly distributed on the Stiefel manifold of $n\times p$ matrices with orthonormal columns.  Then, from \cite[Theorem 2.2.2]{chikuse2012statistics}, we have $\boldsymbol{Q}_k\boldsymbol{Q}_k^\top$ is uniformly distributed on $P_{p,n}$.   Conditional on $\boldsymbol{x}_k$ and all information available immediately before the random directions forming $\boldsymbol{D}_k$ are generated, the result follows from applying Lemma~\ref{lem:generalizedAlphaWellAligned_prob} with $z=p$, $\boldsymbol{x}=\boldsymbol{x}_k$, $\boldsymbol{v}=\nabla f(\boldsymbol{x}_k)$, and $\boldsymbol{X}=\boldsymbol{Q}_k\boldsymbol{Q}_k^\top$.
        
    $\hfill\qed$    
    \end{proof}
    \begin{rema}
        Definition \ref{def:generalizedAlphaWellAligned} is equivalent to requiring the unscaled projection matrix to preserve the fraction $\alpha p\slash n$ of the first-order criticality measure.  Therefore, the factor $p\slash n$ has not disappeared from the convergence analysis (see Lemma \ref{lem:bdpim} below).  Theorem \ref{thm:gen_alphawellaligned} instead shows that a fixed subspace dimension $p$ gives this normalized quality with a fixed probability independent of $n$ and $\boldsymbol{x}_k$.
    \end{rema}
    \begin{rema}
        The probability bound in Theorem \ref{thm:gen_alphawellaligned} may also be written in terms of a failure-probability parameter.  For any integer $p\in(0,n]$, if $\alpha,\beta\in(0,1)$ satisfy
        \begin{equation*}
            \beta>\frac{p+4}{2(p+2)}
            \quad\text{and}\quad
            0<\alpha<1-\sqrt{\frac{2(p+2)(1-\beta)}{p}},
        \end{equation*}
        then Theorem \ref{thm:gen_alphawellaligned} guarantees that $\boldsymbol{D}_k$ is $\alpha$-well-aligned with probability at least $1-\beta$.  Indeed, these conditions imply
        \begin{equation*}
            r_{\alpha,p}=\frac{p(1-\alpha)^2}{2(p+2)}>1-\beta.
        \end{equation*}
        The lower bound on $\beta$ ensures that the stated interval for $\alpha$ is nonempty and reflects the fact that the probability lower bound in Theorem \ref{thm:gen_alphawellaligned} is always less than $1\slash2$.
    \end{rema}

\section{Convergence analysis}\label{sec:converg}
    In this section, we present the convergence analysis of $\myalg$.  The general idea of our convergence analysis is inspired by \cite[Section 3]{chen2024qfully}, which itself is inspired by the general framework proposed in \cite{cartis2023scalable}.  We note that both \cite{cartis2023scalable} and \cite{chen2024qfully} only consider algorithms designed for unconstrained problems.

    It is worth mentioning that all results in this section are independent of the modeling technique used in the algorithm.  In particular, the convergence analysis only requires the models constructed by the algorithm to be $(\mathcal{C},\boldsymbol{Q})$-fully linear with constants $\kappa_{ef}$ and $\kappa_{eg}$ independent of~$k$.  Therefore, although the techniques presented in Section \ref{sec:cqfully} construct linear models, we present the convergence results in terms of quadratic models, which naturally cover linear models by setting all second-order terms to zero. This allows the results in this section to apply in more general settings.  We refer the reader to \cite{audet2017derivative,conn2009introduction} for a detailed discussion on modeling techniques.
    
    Inspired by \cite{cartis2023scalable}, we first classify the iterations into different groups.  Suppose that we have run $\myalg$ for $(K+1)$ iterations and the stopping condition is not triggered, i.e., $\Delta_k\ge\Delta_{\min}$ for all $k\in\{0,\ldots,K\}$.  For $0\le K_1\le K_2\le K$, we denote the following subsets of $\{0,\ldots,K\}$:
    \begin{itemize}
        \item $\mathcal{A}_{[K_1,K_2]}:$ the subset of $\{K_1,\ldots,K_2\}$ where $\boldsymbol{D}_k$ is $\alpha$-well-aligned for $f$ and $\mathcal{C}$ at $\boldsymbol{x}_k$.
        \item $\mathcal{M}_{[K_1,K_2]}:$ the subset of $\{K_1,\ldots,K_2\}$ where the model accuracy test condition is satisfied.
        \item $\mathcal{S}_{[K_1,K_2]}:$ the subset of $\{K_1,\ldots,K_2\}$ where the model accuracy test condition is satisfied with $\rho_k\ge\eta_1$.
        \item $\mathcal{R}_{[K_1,K_2]}:$ the intersection $\mathcal{R}\cap\{K_1,\ldots,K_2\}$ where $\mathcal{R}\subseteq\mathbb{N}$ is the deterministic set of refresh iterations specified in Assumption \ref{ass:subspaceQual}.
    \end{itemize}

    For the remainder of this section, we make use of the following assumptions.
    \begin{assump}\label{ass:smoothness}
        The function $f\in C^{1+}$ with Lipschitz constant $L_{\nabla f}$.
    \end{assump}
    \begin{assump}\label{ass:compactlevelset}
        The function $f$ has a compact feasible level set 
        \begin{equation*}
            \mathcal{C}\cap\mathrm{lev}_{f(\boldsymbol{x}_0)}=\{\boldsymbol{x}:\boldsymbol{x}\in \mathcal{C}~~~\text{and}~~~f(\boldsymbol{x})\le f(\boldsymbol{x}_0)\}.
        \end{equation*}
    \end{assump}
    \begin{assump}\label{ass:bdmodelHess}
        The model Hessian $\|\nabla^2\widehat{m}_k\|\le\kappa_H$ for all $k$, where $\kappa_H\ge 0$ is independent of $k$.
    \end{assump}
    \begin{assump}\label{ass:modelacc}
        Each model $\widehat{m}_k$ belongs to a class of $(\mathcal{C},\boldsymbol{Q}_k)$-fully linear models of $f$ at $\boldsymbol{x}_k$ parameterized by $\Delta_k\in(0,\Delta_{\max}]$ with constants $\kappa_{ef}$ and $\kappa_{eg}$ independent of $k$.
    \end{assump}
    \begin{assump}\label{ass:CauchyDec}
        There exists $\kappa_{\mathrm{tr}}\in (0, 1)$ independent of $k$ such that all solutions $\widehat{\boldsymbol{s}}_k$ of the trust-region subproblem satisfy
        \begin{equation*}
            \widehat{m}_k(\boldsymbol{0}_p)-\widehat{m}_k(\widehat{\boldsymbol{s}}_k)\ge \kappa_{\mathrm{tr}}\pi^m(\boldsymbol{x}_k)\min\left(\frac{\pi^m(\boldsymbol{x}_k)}{\|\nabla^2\widehat{m}_k\|+1}, \Delta_k, 1\right).
        \end{equation*}
    \end{assump}
    \begin{assump}\label{ass:stopInc}
        There exists a constant $\overline{k}\ge 0$ such that $\gamma_{\mathrm{inc}}^k=1$ for all $k\ge\overline{k}$, i.e., we do not increase the trust-region radius when $k\ge\overline{k}$.
    \end{assump}
    \begin{assump}\label{ass:subspaceQual}
        The deterministic refresh schedule $\mathcal{R}\subseteq\mathbb{N}$ supplied to Algorithm~\ref{alg:mainAlg} satisfies the following condition: there exists an integer $T\ge 1$ such that every block of $T$ consecutive iterations contains at least one index in $\mathcal{R}$. At each $k\in\mathcal{R}$, immediately before constructing $\widehat{m}_k$, we let $\boldsymbol{D}_k^U$ be empty and let $\boldsymbol{D}_k=\boldsymbol{D}_k^R$ consist of $p$ random mutually orthogonal directions generated by Algorithm~\ref{alg:dirngen} with $\Delta=\Delta_k$.
    \end{assump}
    \begin{rema}
        Assumption \ref{ass:subspaceQual} is the link between the sample-reuse algorithm and the subspace-quality analysis.  At each $k\in\mathcal{R}$, Theorem \ref{thm:gen_alphawellaligned} gives the required conditional probability that $\boldsymbol{D}_k$ is $\alpha$-well-aligned.  No $\alpha$-well-alignment condition is imposed on iterations outside $\mathcal{R}$.  The minimal full-refresh version (setting $\mathcal{R}=\mathbb{N}$) satisfies Assumption \ref{ass:subspaceQual} with $T=1$.
    \end{rema}
    
    We note that Assumptions \ref{ass:smoothness} to \ref{ass:bdmodelHess} are common assumptions for model-based trust-region algorithms (see, e.g., \cite[Chapter 11]{audet2017derivative}, \cite[Chapter 10]{conn2009introduction}, and \cite{more1983computing}).  Assumption \ref{ass:bdmodelHess} is naturally satisfied if using linear models.  Assumption \ref{ass:modelacc} is the analogue, in the convex-constrained setting, of the standard fully linear model assumption.  Assumption \ref{ass:CauchyDec} states a minimum quality of the solution to the trust-region subproblem.  This condition is motivated by the model decrease achieved at the generalized Cauchy point \cite[Equation (12.2.4)]{conn2000trust} and can be satisfied, e.g., by using \cite[Algorithm~12.2.2]{conn2000trust}.  We refer the reader to \cite[Chapter 12.2.1]{conn2000trust} for further details.  Assumptions \ref{ass:stopInc} and \ref{ass:subspaceQual} can be easily satisfied in practice.

    Under Assumption \ref{ass:compactlevelset}, there exists a constant $f_{\mathrm{low}}\in\mathbb{R}$ such that
    \begin{equation*}
        f(\boldsymbol{x}) \ge f_{\mathrm{low}}~~~\text{for all $\boldsymbol{x}\in \mathcal{C}\cap\mathrm{lev}_{f(\boldsymbol{x}_0)}$}.
    \end{equation*}
    Under Assumption \ref{ass:subspaceQual}, we can develop a lower bound for $|\mathcal{R}_{[K_1,K]}|$, as shown in the following lemma.
    \begin{lemm}\label{lem:countR}
        Suppose that Assumption \ref{ass:subspaceQual} holds.  Let $K\ge K_1\ge\overline{k}\ge 0$.  Then,
        \begin{equation*}
            \left|\mathcal{R}_{[K_1,K]}\right|\ge\left\lfloor\frac{K-K_1+1}{T}\right\rfloor
        \end{equation*}
        where $\lfloor\cdot\rfloor$ is the floor function.
    \end{lemm}
    \begin{proof}
        Let $q=\lfloor(K-K_1+1)\slash T\rfloor$.  The first $qT$ indices in $\{K_1,\ldots,K\}$ can be partitioned into $q$ disjoint blocks of $T$ consecutive iterations.  By Assumption \ref{ass:subspaceQual}, each block contains at least one index in $\mathcal{R}$.  Therefore,
        \begin{equation*}
            \left|\mathcal{R}_{[K_1,K]}\right|\ge q=\left\lfloor\frac{K-K_1+1}{T}\right\rfloor.
        \end{equation*}

    $\hfill\qed$
    \end{proof}

    Because the algorithm has access only to $\pi^m(\boldsymbol{x}_k)$ and not to $\pi^f(\boldsymbol{x}_k)$, the convergence of $\myalg$ depends on establishing a suitable relationship between the two.  Specifically, we require that small values of $\pi^m(\boldsymbol{x}_k)$ imply small values of $\pi^f(\boldsymbol{x}_k)$.  Up to a positive constant, Lemma \ref{lem:bdpim} bounds $\pi^m(\boldsymbol{x}_k)$ below by $\pi^f(\boldsymbol{x}_k)$ for $k\in\mathcal{A}_{[0,K]}\cap\mathcal{M}_{[0,K]}$.
    \begin{lemm}\label{lem:bdpim}
        Suppose that Assumptions \ref{ass:smoothness} and \ref{ass:modelacc} hold.  Let $\epsilon>0$.  For all $k\in\mathcal{A}_{[0,K]}\cap\mathcal{M}_{[0,K]}$ where $\pi^f(\boldsymbol{x}_k)\ge\epsilon$, we have 
        \begin{equation*}
            \pi^m(\boldsymbol{x}_k)\ge\epsilon_g(\epsilon), ~~~\text{where}~~~ \epsilon_g(\epsilon)=\frac{\alpha p\epsilon}{n(\kappa_{eg}\mu+1)}>0.
        \end{equation*}
    \end{lemm}
    \begin{proof}
        The model accuracy test condition is satisfied, so we have $\Delta_k\le\mu\pi^m(\boldsymbol{x}_k)$.  Using the triangle inequality (i.e., $|a+b|\le|a|+|b|$ for all $a,b\in\mathbb{R}$), \cite[Lemma 2.2]{hough2022model}, and Assumption \ref{ass:modelacc}, we have that
        \begin{align*}
            &~~\left|\min\limits_{\substack{\boldsymbol{d}\in \mathcal{C}-\boldsymbol{x}_k\\ \|\boldsymbol{d}\|\le 1}}\nabla f(\boldsymbol{x}_k)^\top \boldsymbol{Q}_k\boldsymbol{Q}_k^\top \boldsymbol{d}\right|\\
            &\le \left|\min\limits_{\substack{\boldsymbol{d}\in \boldsymbol{Q}_k^\top(\mathcal{C}-\boldsymbol{x}_k)\\ \|\boldsymbol{d}\|\le 1}}\nabla f(\boldsymbol{x}_k)^\top \boldsymbol{Q}_k \boldsymbol{d}\right|\\
            &\le  \left|\Bigg|\min\limits_{\substack{\boldsymbol{d}\in \boldsymbol{Q}_k^\top(\mathcal{C}-\boldsymbol{x}_k)\\ \|\boldsymbol{d}\|\le 1}}\nabla f(\boldsymbol{x}_k)^\top \boldsymbol{Q}_k \boldsymbol{d}\Bigg|-\Bigg|\min\limits_{\substack{\boldsymbol{d}\in \boldsymbol{Q}_k^\top(\mathcal{C}-\boldsymbol{x}_k)\\ \|\boldsymbol{d}\|\le 1}}\nabla\widehat{m}_k(\boldsymbol{0}_p)^\top \boldsymbol{d}\Bigg|\right| + \pi^m(\boldsymbol{x}_k)\\
            &\le \max\limits_{\substack{\boldsymbol{d}\in \boldsymbol{Q}_k^\top(\mathcal{C}-\boldsymbol{x}_k)\\ \|\boldsymbol{d}\|\le 1}}\left|\left(\boldsymbol{Q}_k^\top\nabla f(\boldsymbol{x}_k)-\nabla\widehat{m}_k(\boldsymbol{0}_p)\right)^\top \boldsymbol{d}\right| + \pi^m(\boldsymbol{x}_k)\\
            &\le \kappa_{eg}\Delta_k + \pi^m(\boldsymbol{x}_k)\\
            &\le \left(\kappa_{eg}\mu+1\right)\pi^m(\boldsymbol{x}_k).
        \end{align*}

        Since $k\in\mathcal{A}_{[0,K]}$, Definition \ref{def:generalizedAlphaWellAligned} gives
        \begin{align*}
            \left|\min\limits_{\substack{\boldsymbol{d}\in \mathcal{C}-\boldsymbol{x}_k\\ \|\boldsymbol{d}\|\le 1}}\nabla f(\boldsymbol{x}_k)^\top \boldsymbol{Q}_k\boldsymbol{Q}_k^\top \boldsymbol{d}\right| =\frac{p}{n}\left|\min\limits_{\substack{\boldsymbol{d}\in \mathcal{C}-\boldsymbol{x}_k\\ \|\boldsymbol{d}\|\le 1}}\nabla f(\boldsymbol{x}_k)^\top \left(\frac{n}{p}\boldsymbol{Q}_k\boldsymbol{Q}_k^\top\right) \boldsymbol{d}\right| \ge\frac{\alpha p}{n}\pi^f(\boldsymbol{x}_k) \ge\frac{\alpha p}{n}\epsilon
        \end{align*}
        and the result follows.

    $\hfill\qed$
    \end{proof}

    As the definition of $\alpha$-well-aligned matrices essentially indicates that the corresponding normalized projected gradient has a certain quality, we hope that, at least with some probability, there are sufficiently many iterations where $\boldsymbol{D}_k$ is $\alpha$-well-aligned.  We shall make use of the following lemma from \cite[Lemma 4.5]{gratton2015direct}.
    \begin{lemm}\label{lem:GRVZlem4.5}
        Let $r\in(0,1]$ and $(Z_i)_{i=0}^\infty$ be a $\{0,1\}$-valued stochastic process such that
        \begin{equation}\label{cond:P(Zk)gep}
            \mathbb{P}\left[Z_0=1\right]\ge r~~~\text{and}~~~\mathbb{P}\left[Z_i=1\middle|Z_0,\ldots,Z_{i-1}\right]\ge r~~\text{for all $i\ge 1$}.
        \end{equation}
        Let $\lambda\in(0,r)$.  Then,
        \begin{equation*}
            \mathbb{P}\left[\sum_{i=0}^{k-1}Z_i\le\lambda k\right] \le \mathrm{exp}\left(-\frac{(r-\lambda)^2}{2r}k\right).
        \end{equation*}
    \end{lemm}
    \begin{proof}
        See the proof of \cite[Lemma 4.5]{gratton2015direct}.  Notice that the proof therein only uses Condition~\eqref{cond:P(Zk)gep}, instead of the original assumption in the statement (``$p$-probabilistic $\kappa$-descent").

        $\hfill\qed$
    \end{proof}
    
    Using Lemma \ref{lem:GRVZlem4.5}, the following lemma guarantees that the probability of having a relatively small number of iterations where $\boldsymbol{D}_k$ is $\alpha$-well-aligned is bounded above.
    \begin{lemm}\label{lem:probcountA}
        Let $0\le K_1 \le K_2\le K$.  For all $\delta\in (0,1)$, we have 
        \begin{equation*}
            \mathbb{P}\left[\left|\mathcal{A}_{[K_1,K_2]}\cap\mathcal{R}_{[K_1,K_2]}\right|\le(1-\delta)r_{\alpha,p}\left|\mathcal{R}_{[K_1,K_2]}\right|\right]\le e^{-\delta^2r_{\alpha,p}|\mathcal{R}_{[K_1,K_2]}|\slash 2},
        \end{equation*}
        where $r_{\alpha,p}$ is defined in Theorem \ref{thm:gen_alphawellaligned}.
    \end{lemm}
    \begin{proof}
        If $\mathcal{R}_{[K_1,K_2]}$ is empty, then the result is obvious.  Otherwise, suppose that its elements, in increasing order, are $t_0<t_1<\cdots<t_{|\mathcal{R}_{[K_1,K_2]}|-1}$.  For $i\in\{0,1,2,\ldots\}$, we define
        \begin{equation*}
            Z_i = \begin{cases}
                \text{the indicator function of the event $\{t_i\in\mathcal{A}_{[K_1,K_2]}\}$}, &\text{if $i\le |\mathcal{R}_{[K_1,K_2]}|-1$}\\
                1, &\text{if $i\ge |\mathcal{R}_{[K_1,K_2]}|$}.
            \end{cases}
        \end{equation*} 
        Let $\mathcal{F}_{t_i}$ denote all information available immediately before the scheduled refresh at iteration $t_i$.  By Algorithm~\ref{alg:mainAlg} and
        Theorem \ref{thm:gen_alphawellaligned}, $\mathbb{P}[Z_i=1\mid\mathcal{F}_{t_i}] \ge r_{\alpha,p}$ and so $\mathbb{P}[Z_i=1\mid Z_0,\ldots,Z_{i-1}]\ge r_{\alpha,p}$.  Notice that
        \begin{equation*}
            \sum_{i=0}^{|\mathcal{R}_{[K_1,K_2]}|-1}Z_i = \left|\mathcal{A}_{[K_1,K_2]}\cap\mathcal{R}_{[K_1,K_2]}\right|.
        \end{equation*}  
        The result follows from applying Lemma \ref{lem:GRVZlem4.5} with $r=r_{\alpha,p}$, $\lambda=(1-\delta)r_{\alpha,p}$, and $k=|\mathcal{R}_{[K_1,K_2]}|$.

    $\hfill\qed$
    \end{proof}

    Now, we start counting the number of iterations in $\mathcal{A}_{[K_1,K]}$, which, combined with Lemma \ref{lem:probcountA}, will help derive the final convergence result.  As shown in the next two lemmas, we break $\mathcal{A}_{[K_1,K]}$ into two subsets $\mathcal{A}_{[K_1,K]}\cap\mathcal{S}_{[K_1,K]}$ and $\mathcal{A}_{[K_1,K]}\backslash\mathcal{S}_{[K_1,K]}$, and count each of the subsets separately.
    \begin{lemm}\label{lem:countAcapS}
        Suppose that Assumptions \ref{ass:smoothness} to \ref{ass:CauchyDec} hold.  Suppose that $K\ge K_1\ge\overline{k}\ge 0$.  Let $\epsilon>0$.  If $\pi^f(\boldsymbol{x}_k)\ge\epsilon$ for all $k\in\{K_1,\ldots,K\}$, then 
        \begin{equation*}
            \left|\mathcal{A}_{[K_1,K]}\cap\mathcal{S}_{[K_1,K]}\right|\le\phi(\epsilon),
        \end{equation*}
        where
        \begin{equation*}
            \phi(\epsilon)=\frac{f(\boldsymbol{x}_0)-f_{\mathrm{low}}}{\eta_1\kappa_{\mathrm{tr}}\epsilon_g(\epsilon)\min\left(\epsilon_g(\epsilon)\slash(\kappa_H+1), \Delta_{\min}, 1\right)},
        \end{equation*}
        and $\epsilon_g(\cdot)$ is defined in Lemma \ref{lem:bdpim}.
    \end{lemm}
    \begin{proof}
        Let $k\in\mathcal{A}_{[K_1,K]}\cap\mathcal{S}_{[K_1,K]}$.  Then, we have $\rho_k\ge\eta_1$ and $k\in\mathcal{A}_{[K_1,K]}\cap\mathcal{M}_{[K_1,K]}$.  According to Lemma \ref{lem:bdpim} and Assumptions \ref{ass:bdmodelHess} and \ref{ass:CauchyDec}, we have that
        \begin{align*}
            f(\boldsymbol{x}_k)-f(\boldsymbol{x}_{k+1})  &\ge \eta_1\left(\widehat{m}_k(\boldsymbol{0}_p)-\widehat{m}_k(\widehat{\boldsymbol{s}}_k)\right)\\
            &\ge \eta_1\kappa_{\mathrm{tr}}\pi^m(\boldsymbol{x}_k)\min\left(\frac{\pi^m(\boldsymbol{x}_k)}{\left\|\nabla^2\widehat{m}_k\right\|+1},\Delta_k,1\right)\\
            &\ge \eta_1\kappa_{\mathrm{tr}}\epsilon_g(\epsilon)\min\left(\frac{\epsilon_g(\epsilon)}{\kappa_H+1}, \Delta_{\min},1\right).
        \end{align*}
        Notice that $\myalg$ is monotone, all iterates are feasible, and $f$ is bounded below by $f_{\mathrm{low}}$.  We obtain
        \begin{align*}
            f(\boldsymbol{x}_0)-f_{\mathrm{low}} &\ge f(\boldsymbol{x}_{K_1})-f_{\mathrm{low}}\\
            &\ge \sum\limits_{k=K_1}^K\left(f(\boldsymbol{x}_k)-f(\boldsymbol{x}_{k+1})\right)\\
            &\ge \sum\limits_{k\in\mathcal{A}_{[K_1,K]}\cap\mathcal{S}_{[K_1,K]}}\left(f(\boldsymbol{x}_k)-f(\boldsymbol{x}_{k+1})\right)\\
            &\ge \left|\mathcal{A}_{[K_1,K]}\cap\mathcal{S}_{[K_1,K]}\right|\eta_1\kappa_{\mathrm{tr}}\epsilon_g(\epsilon)\min\left(\frac{\epsilon_g(\epsilon)}{\kappa_H+1}, \Delta_{\min}, 1\right),
        \end{align*}
        which completes the proof.

    $\hfill\qed$
    \end{proof}

    \begin{lemm}\label{lem:countAslashS}
        Suppose that Assumption \ref{ass:stopInc} holds.  Suppose that $K\ge K_1\ge\overline{k}$.  Then, 
        \begin{equation*}
            \left|\mathcal{A}_{[K_1,K]}\backslash\mathcal{S}_{[K_1,K]}\right|\le \ln(1\slash\gamma_{\mathrm{dec}})^{-1}\ln(\Delta_{\max}\slash\Delta_{\min}).
        \end{equation*}
    \end{lemm}
    \begin{proof}
        Notice that for all $k\in\mathcal{A}_{[K_1,K]}\backslash\mathcal{S}_{[K_1,K]}$, we have $\Delta_{k+1}=\gamma_{\mathrm{dec}}\Delta_k$.  According to Assumption~\ref{ass:stopInc}, for all $k\ge K_1$, we never increase the trust region.  Since the stopping condition is not triggered and $\Delta_{K_1}\le\Delta_{\max}$, we must have
        \begin{equation*}
            \Delta_{\max}\gamma_{\mathrm{dec}}^{|\mathcal{A}_{[K_1,K]}\backslash\mathcal{S}_{[K_1,K]}|} \ge \Delta_{\min},
        \end{equation*}
        and the result follows from taking the natural logarithm of both sides.

    $\hfill\qed$
    \end{proof}

    Combining Lemmas \ref{lem:countAcapS} and \ref{lem:countAslashS}, we find an upper bound for the number of iterations in $\mathcal{A}_{[K_1,K]}$.
    \begin{lemm}\label{lem:countA}
        Suppose that Assumptions \ref{ass:smoothness} to \ref{ass:stopInc} hold.  Suppose that $K\ge K_1\ge\overline{k}$.  Let $\epsilon>0$.  If $\pi^f(\boldsymbol{x}_k)\ge\epsilon$ for all $k\in\{K_1,\ldots,K\}$, then we have 
        \begin{equation*}
            \left|\mathcal{A}_{[K_1,K]}\right|\le\psi(\epsilon),
        \end{equation*}
        where $\psi(\epsilon)=\phi(\epsilon)+\ln(1\slash\gamma_{\mathrm{dec}})^{-1}\ln(\Delta_{\max}\slash\Delta_{\min})$ and $\phi(\cdot)$ is defined in Lemma \ref{lem:countAcapS}.
    \end{lemm}
    \begin{proof}
        Using Lemmas \ref{lem:countAcapS} and \ref{lem:countAslashS}, we obtain
        \begin{equation*}
            \left|\mathcal{A}_{[K_1,K]}\right| = \left|\mathcal{A}_{[K_1,K]}\cap\mathcal{S}_{[K_1,K]}\right|+\left|\mathcal{A}_{[K_1,K]}\backslash\mathcal{S}_{[K_1,K]}\right| \le \phi(\epsilon)+\ln(1\slash\gamma_{\mathrm{dec}})^{-1}\ln(\Delta_{\max}\slash\Delta_{\min}).
        \end{equation*}

    $\hfill\qed$
    \end{proof}

    Finally, we are ready to present the first convergence result.  The following theorem gives a lower bound on the probability of having an iterate $\boldsymbol{x}_k$ with sufficiently small $\pi^f(\boldsymbol{x}_k)$ in a finite interval of iterations.
    \begin{theo}\label{thm:pim_converg_finite}
        Suppose that Assumptions \ref{ass:smoothness} to \ref{ass:subspaceQual} hold.  Let $\epsilon>0$, $K_1\ge\overline{k}$, and 
        \begin{equation*}
            K\ge \max\left\{K_1, K_1+\frac{2T\psi(\epsilon)}{r_{\alpha,p}}+T-1\right\},
        \end{equation*}
        where $r_{\alpha,p}$ is defined in Theorem \ref{thm:gen_alphawellaligned} and $\psi(\cdot)$ is defined in Lemma \ref{lem:countA}.  If the stopping condition is not triggered for all $k\in\{0,\ldots,K\}$, then we have 
        \begin{equation*}
            \mathbb{P}\left[\min\limits_{K_1\le k\le K}\pi^f(\boldsymbol{x}_k)<\epsilon\right]\ge 1-e^{-r_{\alpha,p}\lfloor(K-K_1+1)\slash T\rfloor\slash 8}.
        \end{equation*}
    \end{theo}
    \begin{proof}
        The lower bound on $K$ implies
        \begin{equation*}
            \psi(\epsilon)\le\frac{1}{2}r_{\alpha,p}\left(\frac{K-K_1+1}{T}-1\right)\le\frac{1}{2}r_{\alpha,p}\left\lfloor\frac{K-K_1+1}{T}\right\rfloor\le\frac{1}{2}r_{\alpha,p}\left|\mathcal{R}_{[K_1,K]}\right|,
        \end{equation*}
        where the last inequality follows from Lemma \ref{lem:countR}.  If $\pi^f(\boldsymbol{x}_k)\ge\epsilon$ for all $k\in\{K_1,\ldots,K\}$, then Lemma \ref{lem:countA} gives
        \begin{equation*}
            \left|\mathcal{A}_{[K_1,K]}\cap\mathcal{R}_{[K_1,K]}\right|\le\left|\mathcal{A}_{[K_1,K]}\right|\le\psi(\epsilon)\le\frac{1}{2}r_{\alpha,p}\left|\mathcal{R}_{[K_1,K]}\right|.
        \end{equation*}
        Therefore, taking $\delta=1\slash2$ in Lemma \ref{lem:probcountA}, we obtain
        \begin{align*}
            \mathbb{P}\left[\min\limits_{K_1\le k\le K}\pi^f(\boldsymbol{x}_k)\ge\epsilon\right]
            &\le\mathbb{P}\left[\left|\mathcal{A}_{[K_1,K]}\cap\mathcal{R}_{[K_1,K]}\right|\le\frac{1}{2}r_{\alpha,p}\left|\mathcal{R}_{[K_1,K]}\right|\right]\\
            &\le e^{-r_{\alpha,p}|\mathcal{R}_{[K_1,K]}|\slash8}\\
            &\le e^{-r_{\alpha,p}\lfloor(K-K_1+1)\slash T\rfloor\slash8}.
        \end{align*}
        The result follows by taking the complement of this event.

    $\hfill\qed$
    \end{proof}

    An almost-sure liminf convergence result can be derived using Theorem \ref{thm:pim_converg_finite}.
    \begin{theo}\label{thm:pim_converg_liminf}
        Suppose that Assumptions \ref{ass:smoothness} to \ref{ass:subspaceQual} hold.  If $\myalg$ is run for infinitely many iterations, then 
        \begin{equation*}
            \mathbb{P}\left[\liminf\limits_{k\to\infty}\pi^f(\boldsymbol{x}_k)=0\right]=1.
        \end{equation*}
    \end{theo}
    \begin{proof}
        Fix an integer $K_1\ge\overline{k}$ and $\epsilon>0$.  Taking $K\to\infty$ in Theorem \ref{thm:pim_converg_finite}, we obtain
        \begin{equation*}
            \mathbb{P}\left[\inf\limits_{k\ge K_1}\pi^f(\boldsymbol{x}_k)<\epsilon\right]=1.
        \end{equation*}
        Taking the countable intersection of these probability-one events over all integers $K_1\ge\overline{k}$ and all $\epsilon=1\slash j$, $j\in\{1,2,\ldots\}$, we obtain that, with probability one, for every $K_1\ge\overline{k}$ and every $j\in\{1,2,\ldots\}$ there exists $k\ge K_1$ such that $\pi^f(\boldsymbol{x}_k)<1\slash j$.  This is equivalent to
        \begin{equation*}
            \mathbb{P}\left[\liminf\limits_{k\to\infty}\pi^f(\boldsymbol{x}_k)=0\right]=1.
        \end{equation*}

    $\hfill\qed$
    \end{proof}

    Finally, we present the complexity result of $\myalg$, i.e., the expected number of iterations required to reach an iterate $\boldsymbol{x}_k$ with $\pi^f(\boldsymbol{x}_k)<\epsilon$.
    \begin{theo}\label{thm:converg_expectation}
        Suppose that Assumptions \ref{ass:smoothness} to \ref{ass:subspaceQual} hold and $\myalg$ is run for infinitely many iterations.  Let $\epsilon>0$ and $K_\epsilon=\min\{k:\pi^f(\boldsymbol{x}_k)<\epsilon\}$.  Denote 
        \begin{equation*}
            K_{\min}(\epsilon) = \left\lceil\max\left\{\overline{k},\overline{k}+\frac{2T\psi(\epsilon)}{r_{\alpha,p}}+T-1\right\}\right\rceil,
        \end{equation*}
        where $r_{\alpha,p}$ is defined in Theorem \ref{thm:gen_alphawellaligned} and $\psi(\cdot)$ is defined in Lemma \ref{lem:countA}.  Then, we have 
        \begin{equation*}
            \mathbb{E}\left[K_\epsilon\right]\le K_{\min}(\epsilon) + \frac{T e^{-r_{\alpha,p}\lfloor(K_{\min}(\epsilon)-\overline{k}+1)\slash T\rfloor\slash 8}}{1-e^{-r_{\alpha,p}\slash 8}} = \mathcal{O}(\epsilon^{-2}).
        \end{equation*}
    \end{theo}
    \begin{proof}
        Taking $K_1=\overline{k}$ in Theorem \ref{thm:pim_converg_finite}, we have, for all $K\ge K_{\min}(\epsilon)$,
        \begin{equation*}
            \mathbb{P}\left[K_\epsilon > K\right] = 1-\mathbb{P}\left[K_\epsilon\le K\right] \le e^{-r_{\alpha,p}\lfloor(K-\overline{k}+1)\slash T\rfloor\slash 8}.
        \end{equation*}
        Using the fact that $\mathbb{E}[X]=\sum_{k=0}^\infty\mathbb{P}[X > k]$ for all non-negative integer-valued random variables~$X$ (e.g., \cite[Problem 3.11.13(a)]{grimmett2020probability}), we obtain
        \begin{align*}
            \mathbb{E}\left[K_\epsilon\right]&=\sum\limits_{K=0}^{K_{\min}(\epsilon)-1}\mathbb{P}\left[K_\epsilon > K\right] + \sum\limits_{K=K_{\min}(\epsilon)}^\infty\mathbb{P}\left[K_\epsilon > K\right]\\
            &\le K_{\min}(\epsilon) + \sum\limits_{K=K_{\min}(\epsilon)}^\infty e^{-r_{\alpha,p}\lfloor(K-\overline{k}+1)\slash T\rfloor\slash 8}\\
            &\le K_{\min}(\epsilon) + T\sum\limits_{j=\lfloor(K_{\min}(\epsilon)-\overline{k}+1)\slash T\rfloor}^\infty e^{-r_{\alpha,p}j\slash8}\\
            &= K_{\min}(\epsilon) + \frac{T e^{-r_{\alpha,p}\lfloor(K_{\min}(\epsilon)-\overline{k}+1)\slash T\rfloor\slash 8}}{1-e^{-r_{\alpha,p}\slash 8}}.
        \end{align*}
        The proof is complete by noticing that $\epsilon_g(\epsilon)$ is linear in $\epsilon$, so $\psi(\epsilon)=\mathcal{O}(\epsilon^{-2})$, while $r_{\alpha,p}>0$ is independent of $\epsilon$.  Hence, $K_{\min}(\epsilon)=\mathcal{O}(\epsilon^{-2})$.

    $\hfill\qed$
    \end{proof}
    
    By noting that the minimal full-refresh version of $\myalg$ must satisfy Assumption \ref{ass:subspaceQual}, we get the following corollary.
    \begin{coroll}\label{cor:minimalCLARSTA}
        Suppose that Assumptions \ref{ass:smoothness} to \ref{ass:stopInc} hold and that the minimal full-refresh version of $\myalg$ is run for infinitely many iterations.  Then, the conclusions of Theorems \ref{thm:pim_converg_liminf} and \ref{thm:converg_expectation} hold with $T=1$.
    \end{coroll}
    \begin{rema}
        Theorem~\ref{thm:converg_expectation} has the same
        $\mathcal{O}(\epsilon^{-2})$ dependence as existing unconstrained
        random subspace trust-region methods
        \cite{cartis2023scalable,chen2024qfully,dzahini2024stochastic}.
        Since Lemma~\ref{lem:bdpim} introduces the factor
        $\alpha p\slash n$, suppressing constants independent of $n$, $p$,
        and $\epsilon$, as well as additional dependence on $p$ in the
        model-accuracy constants, gives
        $\mathcal{O}((n\slash p)^2\epsilon^{-2})$ expected iterations.
        Thus, although $p$ may be chosen independently of $n$, keeping this
        factor uniformly bounded requires $p\slash n$ to remain bounded away
        from zero. For comparison, the stochastic subspace descent method
        in \cite{kozak2021stochastic} achieves
        $\mathcal{O}((n\slash p)\epsilon^{-2})$ iterations in expectation to reach an iterate with $\|\nabla f\|\le\epsilon$, using exact
        directional derivatives, whereas the random subspace model-based DFO
        method in \cite{scheinberg2026complexity} achieves
        $\mathcal{O}(n\epsilon^{-2})$ expected iterations using
        finite-difference models. The original RSDFO high-probability iteration bound $\mathcal{O}(Q_{\max}^4p^2\epsilon^{-2})$ in~\cite{cartis2023scalable} becomes $\mathcal{O}(n^2\epsilon^{-2})$ for the unscaled Johnson-Lindenstrauss transforms considered therein, agreeing with our explicit dimension dependence when $p$ is fixed.  More recently, \cite{cartis2026note} showed that a $1/\sqrt{n}$ rescaling of Johnson-Lindenstrauss transforms in the same RSDFO framework improves this to $\mathcal{O}(n\epsilon^{-2})$ iterations with high probability for $p=\mathcal{O}(1)$ independent of $n$. STARS
        \cite{dzahini2024stochastic} also has an
        $\mathcal{O}(\epsilon^{-2})$ expected iteration bound with $p$
        independent of $n$, although its dimension dependence is contained
        in the subspace and model constants. These results are not directly
        comparable because they use different assumptions and oracle
        information, but they suggest that improving the quadratic dependence
        on $n\slash p$ is an important direction for future work.
    \end{rema}

\section{Numerical experiments}\label{sec:numexps}
    In this section, we design numerical experiments to study the practical performance of $\myalg$ in high dimensions.  For comparison, we also run the numerical tests on the derivative-free model-based trust-region method $\mathrm{COBYLA}$ \cite{Powell1994}.  The implementation of $\mathrm{COBYLA}$ is from the package PDFO \cite{ragonneau2024pdfo}.  We implement $\myalg$ in Python 3, with the random library {\tt numpy.random}.  Some key parameters used in the experiments are set to the following: $\Delta_{\min}=10^{-20}$, $\Delta_{\max}=10^{10}$, $\Delta_0=0.1\max\{\|\boldsymbol{x}_0\|_{\infty},1\}$, $\gamma_{\mathrm{dec}}=0.5$, $\gamma_{\mathrm{inc}}^k=2$ for all $k\in\mathbb{N}$, $\eta_1=0.1$, $\eta_2=0.7$, $\mu=10^{12}$, $\epsilon_{\mathrm{rad}}=10$, and $\epsilon_{\mathrm{geo}}=10^{-6}$.  All parameter values were fixed before the experiments and were used unchanged across all problems and dimensions; no problem-specific tuning was performed.
    
    Since the projection onto $\boldsymbol{Q}_k^\top(\mathcal{C}-\boldsymbol{x}_k)$ is not easy to compute in general, we make the following approximations.
    \begin{itemize}[label=$\bullet$]
        \item We use
            \begin{equation*}
                \widetilde{\pi}^m(\boldsymbol{x}_k)=\left|\min\limits_{\substack{\boldsymbol{d}\in \mathcal{C}-\boldsymbol{x}_k\\ \|\boldsymbol{d}\|\le 1}}\nabla \widehat{m}_k(\boldsymbol{0}_p)^\top \boldsymbol{Q}_k^\top \boldsymbol{d}\right|
            \end{equation*}
            as an approximation of $\pi^m(\boldsymbol{x}_k)$.  Inspired by \cite[Theorem 12.1.4]{conn2000trust}, we approximately compute $\widetilde{\pi}^m(\boldsymbol{x}_k)$ and hence $\pi^m(\boldsymbol{x}_k)$ by
            \begin{equation*}
                \pi^m(\boldsymbol{x}_k) \approx \widetilde{\pi}^m(\boldsymbol{x}_k) \approx \left|\nabla \widehat{m}_k(\boldsymbol{0}_p)^\top \boldsymbol{Q}_k^\top \left(\mathrm{proj}_{\mathcal{C}}\left(\boldsymbol{x}_k-\frac{\boldsymbol{Q}_k\nabla\widehat{m}_k(\boldsymbol{0}_p)}{\|\boldsymbol{Q}_k\nabla\widehat{m}_k(\boldsymbol{0}_p)\|}\right)-\boldsymbol{x}_k\right)\right|.
            \end{equation*}
        \item For the trust-region subproblem, we use the solver provided by \cite{roberts2025model}, which solves the subproblems using the projected gradient descent algorithm \cite[Chapter 10]{beck2017first} with projections onto the two constraint sets computed by Dykstra's projection algorithm \cite{boyle1986method}.  In our implementation, the trust-region subproblem is solved on $\boldsymbol{Q}_k^\top((\mathcal{C}-\boldsymbol{x}_k)\cap\mathrm{col}(\boldsymbol{Q}_k))$ instead of $\boldsymbol{Q}_k^\top(\mathcal{C}-\boldsymbol{x}_k)$.
    \end{itemize}
    Therefore, the experiments are for a practical variant of $\myalg$, not for the generic algorithm analyzed in Section \ref{sec:converg}.  The code of this variant of $\myalg$ is available on GitHub\footnote{https://github.com/yiwchen233/CLARSTA}.  
    
    The purpose of this section is to examine the efficiency of $\myalg$ in high dimensions.  We consider the following objective functions $f:\mathbb{R}^n\to\mathbb{R}$ where $n\in\{10,100,300,600,900,1000\}$ and $x_{(i)}$ denotes the $i$-th component of $\boldsymbol{x}$:
    \begin{itemize}[label=$\bullet$]
        \item ChainRosenbrock: $f(\boldsymbol{x})=\sum\limits_{i=1}^{n-1}\left(100(x_{(i+1)}-x_{(i)}^2)^2+(1-x_{(i)})^2\right)$,
        \item Trigonometric: $f(\boldsymbol{x})=\sum\limits_{i=1}^{n}\left(n-\sum\limits_{j=1}^n\cos x_{(j)}+i\left(1-\cos x_{(i)}\right)-\sin x_{(i)}\right)^2$,
    \end{itemize}
    and the following constraints:
    \begin{itemize}[label=$\bullet$]
        \item ChainRosenbrock:
            \begin{itemize}
                \item box: $\mathcal{C}=[-1,1]^n$,
                \item ball: $\mathcal{C}=B(\boldsymbol{0}_n;\sqrt{n})$,
                \item half-space: $\mathcal{C}=\{\boldsymbol{x}:\boldsymbol{1}_n^\top \boldsymbol{x}\ge 0\}$.
            \end{itemize}
        \item Trigonometric:
            \begin{itemize}
                \item box: $\mathcal{C}=[0,2]^n$,
                \item ball: $\mathcal{C}=B(\boldsymbol{1}_n;\sqrt{n})$,
                \item half-space: $\mathcal{C}=\{\boldsymbol{x}:\boldsymbol{1}_n^\top \boldsymbol{x}\le n\}$.
            \end{itemize}
    \end{itemize}
    Notice that the ChainRosenbrock function has a global minimum $f(\boldsymbol{1}_n)=0$ and the Trigonometric function has a global minimum $f(\boldsymbol{0}_n)=0$.  Therefore, in our experiments, the global minimizer is always at a corner of the box constraint, on the boundary of the ball constraint, and in the interior of the half-space constraint.

    We use $\boldsymbol{x}_0=\boldsymbol{0}_n\in \mathcal{C}$ for the ChainRosenbrock and $\boldsymbol{x}_0=\boldsymbol{1}_n\in \mathcal{C}$ for the Trigonometric as starting points.  We first run $\myalg$ with $p=3$ and $p_{\mathrm{rand}}=3$ until $100(n+1)$ function evaluations are required and record the final objective function value $f^*$.  Then, we run $\mathrm{COBYLA}$ and record the runtime and number of function evaluations required to reach $f^*$.  

    Table \ref{tab:n10-1000_p3prand3} presents the results for $n=10,100,1000$.  In Table \ref{tab:n10-1000_p3prand3}, we use C.R. and Trig. as abbreviations of the ChainRosenbrock and the Trigonometric function, respectively.  We use~$\mathrm{nf}$ to denote the number of function evaluations each algorithm requires to reach $f^*$.  All times are presented in seconds.  We use an asterisk next to the number to indicate that $\mathrm{COBYLA}$ does not reach $f^*$ within $100(n+1)$ function evaluations or $10^5$ seconds, whichever comes first.  If $\mathrm{COBYLA}$ does not reach $f^*$ in $10^5$ seconds, we record its $\mathrm{nf}$ as `N/A', since we have to manually force $\mathrm{COBYLA}$ to stop, and the PDFO package \cite{ragonneau2024pdfo} does not provide the number of function evaluations in this scenario.  
    
    We suppose that the runtime of each algorithm takes the form
    \begin{equation}\label{eq:runtimesplit}
        \mathrm{TotalTime} = \mathrm{AlgTime} + \mathrm{nf}\cdot\mathrm{t}_{\mathrm{feval}},
    \end{equation}
    where $\mathrm{t}_{\mathrm{feval}}$ is the average time for one objective function evaluation and $\mathrm{AlgTime}$ is the runtime required purely by the algorithm (i.e., without function evaluation time).  Since it is not possible to record each objective function evaluation time without influencing $\mathrm{TotalTime}$, we run each algorithm on each test function independently and compute $\mathrm{t}_{\mathrm{feval}}$ by taking the average time of the first $1000$ function evaluations.  In Table \ref{tab:n10-1000_p3prand3}, the $\mathrm{AlgTime}$ is computed using Equation \eqref{eq:runtimesplit}.  

    \begin{table}[htb]
        \caption{Results comparing $\myalg$ ($p=3$, $p_{\mathrm{rand}}=3$) and $\mathrm{COBYLA}$ with $n=10,100,1000$}
        \label{tab:n10-1000_p3prand3}
        \centering
        \begin{tabular}{ccc|ccc|ccc}
        \hline
        \multicolumn{3}{c|}{Problem} & \multicolumn{3}{c|}{$\myalg$} & \multicolumn{3}{c}{$\mathrm{COBYLA}$}      \\
        $n$ & objective & constraint & nf &  AlgTime & TotalTime & nf &  AlgTime & TotalTime \\ \hline
        \multirow{6}{*}{10} & C.R. & box & 1100 & 1.696e+00 & 1.722e+00 & 341 & 5.204e-02 & 5.793e-02 \\
        & C.R. & ball & 1100 & 1.495e+00 & 1.521e+00 & 354 & 7.859e-02 & 8.470e-02 \\
        & C.R. & halfspace & 1100 & 4.377e-01 & 4.637e-01 & 365 & 5.303e-02 & 5.934e-02 \\
        & Trig. & box & 1100 & 1.342e+00 & 1.383e+00 & 1100* & 6.369e-02 & 1.036e-01 \\
        & Trig. & ball & 1100 & 8.489e-01 & 8.895e-01 & 326 & 8.934e-02 & 1.012e-01\\
        & Trig. & halfspace & 1100 & 4.326e-01 & 4.732e-01 & 425 & 5.061e-02 & 6.601e-02 \\ \hline
        \multirow{6}{*}{100} & C.R. & box & 10100 & 1.109e+01 & 1.136e+01 & 4053 & 1.669e+01 & 1.677e+01 \\
        & C.R. & ball & 10100 & 1.097e+01 & 1.124e+01 & 4139 & 3.575e+00 & 3.660e+00 \\
        & C.R. & halfspace & 10100 & 1.114e+01 & 1.140e+01 & 4344 & 3.518e+00 & 3.607e+00 \\
        & Trig. & box & 10100 & 4.014e+01 & 4.959e+01 & 5644 & 3.342e+01 & 3.842e+01 \\
        & Trig. & ball & 10100 & 3.592e+01 & 4.537e+01 & 3506 & 3.277e+00 & 6.382e+00 \\
        & Trig. & halfspace & 10100 & 2.626e+01 & 3.571e+01 & 2494 & 2.060e+00 & 4.269e+00 \\ \hline
        \multirow{6}{*}{1000} & C.R. & box & 100100 & 2.509e+03 & 2.520e+03 & 35673 & 8.367e+04 & 8.367e+04 \\
        & C.R. & ball & 100100 & 2.347e+03 & 2.358e+03 & 36202 & 2.845e+04 & 2.845e+04 \\
        & C.R. & halfspace & 100100 & 2.370e+03 & 2.381e+03 & 34330 & 2.709e+04 & 2.709e+04 \\
        & Trig. & box & 100100 & 3.435e+03 & 1.027e+04 & N/A & N/A & 1.000e+05*\\
        & Trig. & ball & 100100 & 3.405e+03 & 1.024e+04 & 100100* & 7.979e+04 & 8.691e+04 \\
        & Trig. & halfspace & 100100 & 3.456e+03 & 1.029e+04 & 100100* & 7.962e+04 & 8.673e+04 \\ \hline
        \end{tabular}
    \end{table}\FloatBarrier

    As shown in Table \ref{tab:n10-1000_p3prand3}, in low dimensions, $\mathrm{COBYLA}$ is faster than $\myalg$ and also requires fewer function evaluations to reach $f^*$.  This is not surprising, as $\mathrm{COBYLA}$ constructs full-space linear interpolation models to approximate both the objective function and constraints, while $\myalg$ only constructs linear approximation models in a $p$-dimensional subspace.  In low dimensions, full-space models are computationally cheap to construct, so the performance of $\mathrm{COBYLA}$ is better than $\myalg$.
    
    In high dimensions, however, $\myalg$ is much faster than $\mathrm{COBYLA}$ to reach $f^*$.  Indeed, Table \ref{tab:n10-1000_p3prand3} shows that the $\mathrm{AlgTime}$ of $\mathrm{COBYLA}$ increases much more rapidly than that of $\myalg$ as $n$ increases.
    
    It should be noted that, however, despite the rather low efficiency in high dimensions, $\mathrm{COBYLA}$ still has the potential to take fewer function evaluations than $\myalg$.  This implies that, in high dimensions, $\myalg$ is likely to be more efficient than $\mathrm{COBYLA}$ when the function evaluation is not expensive (i.e., small $\mathrm{t}_{\mathrm{feval}}$), and $\mathrm{COBYLA}$ is likely to be more efficient when the function evaluation is expensive (i.e., large $\mathrm{t}_{\mathrm{feval}}$).  

    In order to explore how $\mathrm{t}_{\mathrm{feval}}$ influences the relative performance of $\myalg$ and $\mathrm{COBYLA}$, we use the data for $n\in\{100,300,600,900,1000\}$ to find the smallest value of $\mathrm{t}^*_{\mathrm{feval}}>0$ such that
    \begin{equation*}
        \mathrm{AlgTime}_{\myalg} + \mathrm{nf}_{\myalg}\mathrm{t}^*_{\mathrm{feval}} = \mathrm{AlgTime}_{\mathrm{COBYLA}} + \mathrm{nf}_{\mathrm{COBYLA}}\mathrm{t}^*_{\mathrm{feval}}.
    \end{equation*}
    Notice that in our experiments, $\mathrm{nf}_{\mathrm{COBYLA}}\le\mathrm{nf}_{\myalg}=100(n+1)$.  We list below the few possible cases and how we present the result for each case in Table \ref{tab:fevaltime_p3prand3}.  All numerical values are presented in seconds.
    \begin{itemize}[label=$\bullet$]
        \item $\mathrm{AlgTime}_{\myalg}<\mathrm{AlgTime}_{\mathrm{COBYLA}}$ and $\mathrm{nf}_{\mathrm{COBYLA}}<\mathrm{nf}_{\myalg}$.\\
            In this case, $\mathrm{t}^*_{\mathrm{feval}}>0$ exists, $\myalg$ is more efficient than $\mathrm{COBYLA}$ when $\mathrm{t}_{\mathrm{feval}}<\mathrm{t}^*_{\mathrm{feval}}$, and $\mathrm{COBYLA}$ is more efficient than $\myalg$ when $\mathrm{t}_{\mathrm{feval}}>\mathrm{t}^*_{\mathrm{feval}}$.  In Table \ref{tab:fevaltime_p3prand3}, we record this $\mathrm{t}^*_{\mathrm{feval}}>0$.
        \item $\mathrm{AlgTime}_{\myalg}\ge\mathrm{AlgTime}_{\mathrm{COBYLA}}$ and $\mathrm{nf}_{\mathrm{COBYLA}}<\mathrm{nf}_{\myalg}$.\\
            In this case, $\mathrm{t}^*_{\mathrm{feval}}>0$ does not exist.  In fact, now $\myalg$ is no faster and requires more function evaluations than $\mathrm{COBYLA}$, which implies that $\mathrm{COBYLA}$ is more efficient than $\myalg$, regardless of $\mathrm{t}_{\mathrm{feval}}$.  In Table \ref{tab:fevaltime_p3prand3}, we record this case as `0'.
        \item $\mathrm{nf}_{\mathrm{COBYLA}}$ is recorded as `N/A' or $\mathrm{nf}_{\mathrm{COBYLA}}=\mathrm{nf}_{\myalg}=100(n+1)$.\\
            In this case, $\mathrm{COBYLA}$ does not successfully find a solution that is as good as $f^*$ within $100(n+1)$ function evaluations or $10^5$ seconds.  Since $\mathrm{TotalTime}_{\mathrm{COBYLA}} > \mathrm{TotalTime}_{\myalg}$ in all such cases, this means $\myalg$ is more efficient regardless of $\mathrm{t}_{\mathrm{feval}}$.  Therefore, we record this case as `+Inf' in Table \ref{tab:fevaltime_p3prand3}.
    \end{itemize}

    \begin{table}[htb]
        \caption{Maximum $\mathrm{t}_{\mathrm{feval}}$ for $\myalg$ ($p=3$, $p_{\mathrm{rand}}=3$) to be more efficient than $\mathrm{COBYLA}$ in high dimensions}
        \label{tab:fevaltime_p3prand3}
        \centering
        \begin{tabular}{ccccccc}
        \hline
        objective & constraint & $n=100$ & $n=300$ & $n=600$ & $n=900$ & $n=1000$\\ \hline
        ChainRosenbrock & box & 9.261e-04 & 4.389e-02 & 4.503e-01 & +Inf & 1.260e+00 \\
        ChainRosenbrock & ball & 0 & 7.622e-03 & 1.277e-01 & 1.373e+00 & 4.085e-01 \\
        ChainRosenbrock & halfspace & 0 & 7.647e-03 & 1.278e-01 & 1.372e+00 & 3.758e-01\\
        Trigonometric & box & 0 & 1.475e-01 & 2.872e+00 & +Inf & +Inf\\
        Trigonometric & ball & 0 & +Inf & 1.717e+01 & +Inf & +Inf\\
        Trigonometric & halfspace & 0 & 0 & +Inf & +Inf & +Inf\\ \hline
        \end{tabular}
    \end{table}\FloatBarrier
    
    Table \ref{tab:fevaltime_p3prand3} clearly illustrates the potential of $\myalg$ in high dimensions.

    For completeness, we repeat the experiments with the parameters $p$ and $p_{\mathrm{rand}}$ in $\myalg$ set to $(p,p_{\mathrm{rand}})=(1,1)$ and $(p,p_{\mathrm{rand}})=(3,1)$.  The corresponding results are presented in Tables~\ref{tab:n10-10000_p1prand1} and \ref{tab:n10-1000_p3prand1} in Appendix \ref{app:morenumexps}.  Overall, Tables~\ref{tab:n10-10000_p1prand1} and \ref{tab:n10-1000_p3prand1} exhibit trends similar to those observed in Table~\ref{tab:n10-1000_p3prand3}.  Similarly, Tables~\ref{tab:fevaltime_p1prand1} and \ref{tab:fevaltime_p3prand1} show that the range of $\mathrm{t}_{\mathrm{feval}}$ for which $\myalg$ is more efficient than $\mathrm{COBYLA}$ generally increases with the problem dimension, with several cases becoming `+Inf'.  It is worth mentioning that for $(p,p_{\mathrm{rand}})=(1,1)$, we also run the experiment on the ChainRosenbrock function with $n=10000$.  The difference in efficiency is more pronounced when $n=10000$: $\myalg$ takes about half an hour to reach $f^*$, whereas $\mathrm{COBYLA}$ cannot reach $f^*$ within $10^5$ seconds in any of the tests.
    
    Comparing Tables \ref{tab:n10-1000_p3prand3}, \ref{tab:n10-10000_p1prand1}, and \ref{tab:n10-1000_p3prand1}, we see that although $\myalg$ becomes slower when $p$ and/or $p_{\mathrm{rand}}$ are larger than~$1$, it remains significantly faster than $\mathrm{COBYLA}$.  Moreover, the number of function evaluations required by $\mathrm{COBYLA}$ to reach the solution obtained by $\myalg$ increases, indicating that $\myalg$ produces higher-quality solutions when $p$ and/or $p_{\mathrm{rand}}$ are slightly larger than $1$.  These observations are consistent with the numerical results reported in~\cite{chen2024qfully}, which compare $\mathrm{QARSTA}$, an algorithm similar to $\myalg$ but designed for unconstrained problems, under different choices of $p$ and $p_{\mathrm{rand}}$.

\section{Conclusion}\label{sec:conclusion}
    In this paper, we propose $\myalg$, a random subspace trust-region algorithm for general convex-constrained derivative-free optimization problems.  We define a class of $(\mathcal{C},\boldsymbol{Q})$-fully linear models, which is only required to provide the same level of accuracy as fully linear models~\cite{conn2009global} on the projection of the constraint set $\mathcal{C}$ onto the subspace determined by matrix $\boldsymbol{Q}$.  This definition is a generalization of fully linear models \cite{conn2009global}, $\mathcal{C}$-pointwise fully linear models \cite{hough2022model,roberts2025model}, and $\boldsymbol{Q}$-fully linear models~\cite{cartis2023scalable}.  A new geometry measure of sample sets is introduced, and techniques to construct and manage a class of $(\mathcal{C},\boldsymbol{Q})$-fully linear models are provided.

    We define a normalized convex-constrained version of the $\alpha$-well-aligned matrices introduced in \cite{cartis2023scalable}.  Matrices that satisfy this new definition can be used to construct subspaces that preserve the first-order criticality measure for convex-constrained optimization \cite{conn2000trust}.  Using moment properties of uniformly sampled orthogonal projection matrices, we provide a sampling technique for such matrices.  We show that the resulting matrices satisfy the proposed definition with a probability lower bound independent of the ambient dimension.

    Based on these new results, we prove that the first-order criticality measure converges to zero in the liminf almost surely and provide an $\mathcal{O}(\epsilon^{-2})$ expected iteration complexity bound for $\myalg$.  We design numerical experiments to compare the performance of $\myalg$ and $\mathrm{COBYLA}$ \cite{Powell1994} (from the package PDFO \cite{ragonneau2024pdfo}) in both low and high dimensions. Numerical results on problems with dimensions up to 10000 demonstrate the reliable performance of $\myalg$ in high dimensions.

    An obvious next step of this paper is to extend the implementation of $\myalg$ to use quadratic models.  This will require a management procedure for quadratic models.  To this end, the modeling technique introduced in \cite{chen2024qfully} may be used to help.  Another next step is to improve the dimension dependence of the constants in the convergence and complexity bounds.

\section*{\bf Data availability and declarations}
All data are available upon request from the corresponding author.  The authors have no competing interests to declare that are relevant to the content of this article.

\appendix
\section{High-level flowchart of the minimal full-refresh version of $\myalg$}\label{app:flowchart}
    \begin{figure}[htbp!]
    \centering
        \begin{tikzpicture}[
            node distance=6mm,
            >=Latex,
            block/.style={
                rectangle,
                rounded corners,
                draw,
                align=left,
                inner sep=4pt,
                minimum width=30mm
            },
            decision/.style={
                diamond,
                draw,
                align=center,
                aspect=3,
                inner sep=2pt
            },
            term/.style={
                rounded rectangle,
                draw,
                align=center,
                inner sep=4pt
            },
            line/.style={-Latex}
        ]
    
            % Nodes
            \node[term] (start) {\textbf{Start}};
    
            \node[block,below=of start] (init) {
            \textbf{Initialize direction matrix}\\
            Generate $p$ random mutually orthogonal directions\\
            by Alg.~\ref{alg:dirngen} and set
            $\boldsymbol{D}_0=\boldsymbol{D}_0^R$
            };
    
            \node[block,below=of init] (model) {
            \textbf{Iteration $k=0,1,2,\dots$}\\
            Construct subspace model $\widehat{m}_k$
            };
    
            \node[decision,below=of model] (test) {
            \textbf{Model accuracy test}\\
            $\Delta_k\le\mu\pi^m(\boldsymbol{x}_k)$
            };
    
            % Pass branch
            \node[block,below=of test] (trs) {
            \textbf{Trust-region step}\\
            Solve subspace trust-region subproblem\\
            Project step to $\mathcal{C}$
            };
    
            \node[block,below=of trs] (accept) {
            \textbf{Trust-region update}\\
            Calculate ratio $\rho_k$\\
            Update $\Delta_{k+1}$ and $\boldsymbol{x}_{k+1}$
            };
    
            % Fail branch
            \node[block,right=18mm of test] (shrink) {
            \textbf{Shrink step}\\
            $\Delta_{k+1}=\gamma_{\mathrm{dec}}\Delta_k$\\
            $\boldsymbol{x}_{k+1}=\boldsymbol{x}_k$
            };
    
            % Full refresh
            \node[block,below=of accept] (dirs) {
            \textbf{Generate new direction matrix}\\
            Generate $p$ fresh random mutually orthogonal directions\\
            by Alg.~\ref{alg:dirngen} with $\Delta=\Delta_{k+1}$\\
            Set $\boldsymbol{D}_{k+1}=\boldsymbol{D}_{k+1}^R$
            };
    
            % Stop
            \node[decision,below=of dirs] (stoptest) {
            \textbf{Stopping test}\\
            $\Delta_{k+1}<\Delta_{\min}$
            };
    
            \node[term,below=of stoptest] (stop) {\textbf{Stop}};
    
            % Lines
            \draw[line] (start) -- (init);
            \draw[line] (init) -- (model);
            \draw[line] (model) -- (test);
    
            \draw[line] (test) -- node[right]{Pass} (trs);
            \draw[line] (trs) -- (accept);
            \draw[line] (accept) -- (dirs);
    
            \draw[line] (test.east) -- node[above]{Fail} (shrink.west);
            \draw[line] (shrink.south) |- (dirs.east);
    
            \draw[line] (dirs) -- (stoptest);
            \draw[line] (stoptest) -- node[right]{Pass} (stop);
            \draw[line] (stoptest.west)
                -- node[above]{Fail} ++(-30mm,0)
                |- (model.west);
    
        \end{tikzpicture}
    \end{figure}\FloatBarrier

\section{More numerical results}\label{app:morenumexps}
    \begin{table}[htb]
        \caption{Results comparing $\myalg$ ($p=1$, $p_{\mathrm{rand}}=1$) and $\mathrm{COBYLA}$ with $n=10,100,1000,10000$}
        \label{tab:n10-10000_p1prand1}
        \centering
        \begin{tabular}{ccc|ccc|ccc}
        \hline
        \multicolumn{3}{c|}{Problem} & \multicolumn{3}{c|}{$\myalg$} & \multicolumn{3}{c}{$\mathrm{COBYLA}$}      \\
        $n$ & objective & constraint & nf &  AlgTime & TotalTime & nf &  AlgTime & TotalTime \\ \hline
        \multirow{6}{*}{10} & C.R. & box & 1100 & 1.536e+00 & 1.562e+00  & 341  & 8.444e-03 & 1.434e-02  \\
        & C.R. & ball & 1100 & 1.582e+00 & 1.608e+00  & 340  & 2.806e-02 & 3.394e-02  \\
        & C.R. & halfspace & 1100 & 6.591e-01 & 6.850e-01  & 274 & 7.786e-03 & 1.252e-02  \\
        & Trig. & box & 1100 & 1.147e+00 & 1.187e+00  & 1100* & 1.517e-02 & 5.503e-02  \\
        & Trig. & ball & 1100 & 1.220e+00 & 1.261e+00  & 347  & 2.876e-02 & 4.133e-02  \\
        & Trig. & halfspace & 1100 & 1.013e+00 & 1.054e+00  & 165  & 7.900e-03 & 1.389e-02  \\ \hline
        % \end{tabular}
        % \end{table}
    
        % \begin{table}[htb]
        % \caption{Results comparing $\myalg$ and $\mathrm{COBYLA}$ with $n=100$}
        % \label{tab:n100}
        % \centering
        % \begin{tabular}{cc|ccc|ccc}
        % \hline
        % \multicolumn{2}{c|}{Problem} & \multicolumn{3}{c|}{$\myalg$} & % \multicolumn{3}{c}{$\mathrm{COBYLA}$}      \\
        % objective & constraint & nf &  AlgTime & TotalTime & nf &  AlgTime & TotalTime \\ \hline
        \multirow{6}{*}{100} & C.R. & box & 10100 & 5.481e+00 & 5.748e+00 & 3065 & 1.224e+01 & 1.230e+01\\
        & C.R. & ball & 10100 & 5.662e+00 & 5.929e+00 & 3113 & 2.659e+00 & 2.723e+00\\
        & C.R. & halfspace & 10100 & 5.864e+00 & 6.131e+00 & 3121 & 2.475e+00 & 2.540e+00\\
        & Trig. & box & 10100 & 1.408e+01 & 2.353e+01 & 2524 & 1.502e+01 & 1.726e+01\\
        & Trig. & ball & 10100 & 1.538e+01 & 2.483e+01 & 3167 & 2.837e+00 & 5.642e+00\\
        & Trig. & halfspace & 10100 & 7.191e+00 & 1.664e+01 & 2879 & 2.484e+00 & 5.034e+00\\ \hline
        % \end{tabular}
        % \end{table}
    
        % \begin{table}[htb]
        % \caption{Results comparing $\myalg$ and $\mathrm{COBYLA}$ with $n=1000$}
        % \label{tab:n1000}
        % \centering
        % \begin{tabular}{cc|ccc|ccc}
        % \hline
        % \multicolumn{2}{c|}{Problem} & \multicolumn{3}{c|}{$\myalg$} & \multicolumn{3}{c}{$\mathrm{COBYLA}$}      \\
        % objective & constraint & nf &  AlgTime & TotalTime & nf &  AlgTime & TotalTime \\ \hline
        \multirow{6}{*}{1000} & C.R. & box & 100100 & 5.359e+01 & 6.451e+01 & 26716 & 6.196e+04 & 6.196e+04\\
        & C.R. & ball & 100100 & 5.821e+01 & 6.913e+01 & 26337 & 2.048e+04 & 2.048e+04\\
        & C.R. & halfspace & 100100 & 6.257e+01 & 7.349e+01 & 27054 & 2.113e+04 & 2.113e+04\\
        & Trig. & box & 100100 & 1.446e+02 & 6.975e+03 & N/A & N/A & 1.000e+05*\\
        & Trig. & ball & 100100 & 2.219e+02 & 7.052e+03 & 100100* & 7.967e+04 & 8.678e+04\\
        & Trig. & halfspace & 100100 & 7.082e+01 & 6.901e+03 & 100100* & 7.944e+04 & 8.655e+04\\ \hline
        \multirow{3}{*}{10000} & C.R. & box & 1000100 & 1.091e+03 & 1.836e+03 & N/A & N/A & 1.000e+05*\\
        & C.R. & ball & 1000100 & 1.215e+03 & 1.960e+03 & N/A & N/A & 1.000e+05*\\
        & C.R. & halfspace & 1000100 & 1.183e+03 & 1.928e+03 & N/A & N/A & 1.000e+05*\\ \hline
        \end{tabular}
    \end{table}\FloatBarrier

    \begin{table}[htb]
        \caption{Maximum $\mathrm{t}_{\mathrm{feval}}$ for $\myalg$ ($p=1$, $p_{\mathrm{rand}}=1$) to be more efficient than $\mathrm{COBYLA}$ in high dimensions}
        \label{tab:fevaltime_p1prand1}
        \centering
        \begin{tabular}{ccccccc}
        \hline
        objective & constraint & $n=100$ & $n=300$ & $n=600$ & $n=900$ & $n=1000$\\ \hline
        ChainRosenbrock & box & 9.605e-04 & 3.447e-02 & 2.881e-01 & 1.921e+00 & 8.435e-01\\
        ChainRosenbrock & ball & 0 & 8.727e-03 & 9.130e-02 & 6.324e-01 & 2.769e-01\\
        ChainRosenbrock & halfspace & 0 & 7.775e-03 & 9.117e-02 & 6.182e-01 & 2.884e-01\\ 
        Trigonometric & box & 1.241e-04 & 1.339e-01 & 2.526e+00 & 1.375e+03 & +Inf\\
        Trigonometric & ball & 0 & 5.226e-02 & +Inf & +Inf & +Inf\\
        Trigonometric & halfspace & 0 & 3.057e-02 & +Inf & +Inf & +Inf\\ \hline
        \end{tabular}
    \end{table}\FloatBarrier
    
    \begin{table}[htb]
        \caption{Results comparing $\myalg$ ($p=3$, $p_{\mathrm{rand}}=1$) and $\mathrm{COBYLA}$ with $n=10,100,1000$}
        \label{tab:n10-1000_p3prand1}
        \centering
        \begin{tabular}{ccc|ccc|ccc}
        \hline
        \multicolumn{3}{c|}{Problem} & \multicolumn{3}{c|}{$\myalg$} & \multicolumn{3}{c}{$\mathrm{COBYLA}$}      \\
        $n$ & objective & constraint & nf &  AlgTime & TotalTime & nf &  AlgTime & TotalTime \\ \hline
        \multirow{6}{*}{10} & C.R. & box & 1100 & 3.175e+00 & 3.201e+00 & 583 & 5.367e-02 & 6.375e-02 \\
        & C.R. & ball & 1100 & 2.823e+00 & 2.849e+00 & 567 & 9.005e-02 & 9.985e-02 \\
        & C.R. & halfspace & 1100 & 8.618e-01 & 8.877e-01 & 583 & 5.245e-02 & 6.252e-02 \\
        & Trig. & box & 1100 & 2.834e+00 & 2.875e+00 & 213 & 5.222e-02 & 5.994e-02 \\
        & Trig. & ball & 1100 & 1.306e+00 & 1.347e+00 & 1100* & 1.251e-01 & 1.650e-01 \\
        & Trig. & halfspace & 1100 & 7.747e-01 & 8.153e-01 & 1100* & 5.932e-02 & 9.919e-02 \\ \hline
        \multirow{6}{*}{100} & C.R. & box & 10100 & 2.200e+01 & 2.227e+01 & 5542 & 2.339e+01 & 2.350e+01 \\
        & C.R. & ball & 10100 & 2.090e+01 & 2.116e+01 & 5924 & 5.120e+00 & 5.242e+00 \\
        & C.R. & halfspace & 10100 & 2.039e+01 & 2.065e+01 & 5688 & 4.598e+00 & 4.715e+00 \\
        & Trig. & box & 10100 & 5.583e+01 & 6.528e+01 & 5764 & 3.303e+01 & 3.814e+01 \\
        & Trig. & ball & 10100 & 5.050e+01 & 5.995e+01 & 4112 & 3.817e+00 & 7.458e+00 \\
        & Trig. & halfspace & 10100 & 3.438e+01 & 4.383e+01 & 4893 & 4.406e+00 & 8.740e+00 \\ \hline
        \multirow{6}{*}{1000} & C.R. & box & 100100 & 4.666e+03 & 4.677e+03 & N/A & N/A & 1.000e+05*\\
        & C.R. & ball & 100100 & 4.982e+03 & 4.993e+03 & 44488 & 3.513e+04 & 3.514e+04 \\
        & C.R. & halfspace & 100100 & 4.996e+03 & 5.007e+03 & 44342 & 3.506e+04 & 3.506e+04 \\
        & Trig. & box & 100100 & 7.108e+03 & 1.394e+04 & N/A & N/A & 1.000e+05*\\
        & Trig. & ball & 100100 & 6.843e+03 & 1.367e+04 & 100100* & 7.976e+04 & 8.688e+04 \\
        & Trig. & halfspace & 100100 & 6.771e+03 & 1.360e+04 & 100100* & 7.973e+04 & 8.684e+04 \\ \hline
        \end{tabular}
    \end{table}\FloatBarrier

    \begin{table}[htb]
        \caption{Maximum $\mathrm{t}_{\mathrm{feval}}$ for $\myalg$ ($p=3$, $p_{\mathrm{rand}}=1$) to be more efficient than $\mathrm{COBYLA}$ in high dimensions}
        \label{tab:fevaltime_p3prand1}
        \centering
        \begin{tabular}{ccccccc}
        \hline
        objective & constraint & $n=100$ & $n=300$ & $n=600$ & $n=900$ & $n=1000$\\ \hline
        ChainRosenbrock & box & 3.033e-04 & 5.920e-02 & 6.997e-01 & +Inf & +Inf\\
        ChainRosenbrock & ball & 0 & 8.110e-03 & 1.760e-01 & 2.873e+00 & 5.421e-01\\
        ChainRosenbrock & halfspace & 0 & 7.367e-03 & 1.760e-01 & 2.835e+00 & 5.391e-01\\ 
        Trigonometric & box & 0 & 1.518e-01 & +Inf & +Inf & +Inf\\
        Trigonometric & ball & 0 & 0 & +Inf & +Inf & +Inf\\
        Trigonometric & halfspace & 0 & 0 & +Inf & +Inf & +Inf\\ \hline
        \end{tabular}
    \end{table}\FloatBarrier

\bibliographystyle{siam}
\bibliography{references}

\begin{thebibliography}{10}

\bibitem{abramson2006convergence}
{\sc M.~A. Abramson and C.~Audet}, {\em Convergence of mesh adaptive direct search to second-order stationary points}, SIAM Journal on Optimization, 17 (2006), pp.~606--619.

\bibitem{alarie2021two}
{\sc S.~Alarie, C.~Audet, A.~E. Gheribi, M.~Kokkolaras, and S.~Le~Digabel}, {\em Two decades of blackbox optimization applications}, EURO Journal on Computational Optimization, 9 (2021), p.~100011.

\bibitem{audet2006mesh}
{\sc C.~Audet and J.~E. Dennis~Jr}, {\em Mesh adaptive direct search algorithms for constrained optimization}, SIAM Journal on optimization, 17 (2006), pp.~188--217.

\bibitem{audet2008parallel}
{\sc C.~Audet, J.~E. Dennis~Jr, and S.~Le~Digabel}, {\em Parallel space decomposition of the mesh adaptive direct search algorithm}, SIAM Journal on Optimization, 19 (2008), pp.~1150--1170.

\bibitem{audet2017derivative}
{\sc C.~Audet and W.~Hare}, {\em Derivative-free and blackbox optimization}, Springer, Cham, 2017.

\bibitem{beck2017first}
{\sc A.~Beck}, {\em First-order methods in optimization}, SIAM, 2017.

\bibitem{boyle1986method}
{\sc J.~P. Boyle and R.~L. Dykstra}, {\em A method for finding projections onto the intersection of convex sets in {H}ilbert spaces}, in Advances in Order Restricted Statistical Inference, Springer, 1986, pp.~28--47.

\bibitem{cartis2022randomiseda}
{\sc C.~Cartis, J.~Fowkes, and Z.~Shao}, {\em A randomised subspace {G}auss-{N}ewton method for nonlinear least-squares}, 2022.
\newblock \href{https://arxiv.org/abs/2211.05727}{arXiv:2211.05727}.

\bibitem{cartis2022randomisedb}
\leavevmode\vrule height 2pt depth -1.6pt width 23pt, {\em Randomised subspace methods for non-convex optimization, with applications to nonlinear least-squares}, 2022.
\newblock \href{https://arxiv.org/abs/2211.09873}{arXiv:2211.09873}.

\bibitem{cartis2023scalable}
{\sc C.~Cartis and L.~Roberts}, {\em Scalable subspace methods for derivative-free nonlinear least-squares optimization}, Mathematical Programming, 199 (2023), pp.~461--524.

\bibitem{cartis2024randomized}
\leavevmode\vrule height 2pt depth -1.6pt width 23pt, {\em Randomized subspace derivative-free optimization with quadratic models and second-order convergence}, 2024.
\newblock \href{https://arxiv.org/abs/2412.14431}{arXiv:2412.14431}.

\bibitem{cartis2026note}
\leavevmode\vrule height 2pt depth -1.6pt width 23pt, {\em A note on the complexity of random subspace model-based methods for derivative-free optimization}, 2026.
\newblock \href{https://arxiv.org/abs/2608.17307}{arXiv:2608.17307}.

\bibitem{chen2024qfully}
{\sc Y.~Chen, W.~Hare, and A.~Wiebe}, {\em {$Q$}-fully quadratic modeling and its application in a random subspace derivative-free method}, Computational Optimization and Applications, 89 (2024), pp.~317--360.

\bibitem{chikuse2012statistics}
{\sc Y.~Chikuse}, {\em Statistics on special manifolds}, vol.~174, Springer, 2012.

\bibitem{conejo2013global}
{\sc P.~D. Conejo, E.~W. Karas, L.~G. Pedroso, A.~A. Ribeiro, and M.~Sachine}, {\em Global convergence of trust-region algorithms for convex constrained minimization without derivatives}, Applied Mathematics and Computation, 220 (2013), pp.~324--330.

\bibitem{conn2000trust}
{\sc A.~R. Conn, N.~I.~M. Gould, and P.~L. Toint}, {\em Trust region methods}, SIAM, 2000.

\bibitem{conn2008geometry}
{\sc A.~R. Conn, K.~Scheinberg, and L.~N. Vicente}, {\em Geometry of interpolation sets in derivative free optimization}, Mathematical programming, 111 (2008), pp.~141--172.

\bibitem{conn2009global}
\leavevmode\vrule height 2pt depth -1.6pt width 23pt, {\em Global convergence of general derivative-free trust-region algorithms to first-and second-order critical points}, SIAM Journal on Optimization, 20 (2009), pp.~387--415.

\bibitem{conn2009introduction}
\leavevmode\vrule height 2pt depth -1.6pt width 23pt, {\em Introduction to derivative-free optimization}, SIAM, 2009.

\bibitem{dzahini2024revisiting}
{\sc K.~J. Dzahini, F.~Rinaldi, C.~W. Royer, and D.~Zeffiro}, {\em Revisiting theoretical guarantees of direct-search methods}, 2024.
\newblock \href{https://arxiv.org/abs/2403.05322}{arXiv:2403.05322}.

\bibitem{dzahini2024stochastic}
{\sc K.~J. Dzahini and S.~M. Wild}, {\em Stochastic trust-region algorithm in random subspaces with convergence and expected complexity analyses}, SIAM Journal on Optimization, 34 (2024), pp.~2671--2699.

\bibitem{eaton2007multivariate}
{\sc M.~L. Eaton}, {\em Multivariate statistics: a vector space approach}, vol.~53 of Lecture Notes - Monograph Series, Institute of Mathematical Statistics, 2007.

\bibitem{feurer2019hyperparameter}
{\sc M.~Feurer and F.~Hutter}, {\em Hyperparameter optimization}, Springer, Cham, 2019, pp.~3--33.

\bibitem{ghanbari2017black}
{\sc H.~Ghanbari and K.~Scheinberg}, {\em Black-box optimization in machine learning with trust region based derivative free algorithm}, 2017.
\newblock \href{https://arxiv.org/abs/1703.06925}{arXiv:1703.06925}.

\bibitem{gratton2015direct}
{\sc S.~Gratton, C.~W. Royer, L.~N. Vicente, and Z.~Zhang}, {\em Direct search based on probabilistic descent}, SIAM Journal on Optimization, 25 (2015), pp.~1515--1541.

\bibitem{gratton2019direct}
\leavevmode\vrule height 2pt depth -1.6pt width 23pt, {\em Direct search based on probabilistic feasible descent for bound and linearly constrained problems}, Computational Optimization and Applications, 72 (2019), pp.~525--559.

\bibitem{gratton2011active}
{\sc S.~Gratton, P.~L. Toint, and A.~Tr{\"o}ltzsch}, {\em An active-set trust-region method for derivative-free nonlinear bound-constrained optimization}, Optimization Methods and Software, 26 (2011), pp.~873--894.

\bibitem{grimmett2020probability}
{\sc G.~Grimmett and D.~Stirzaker}, {\em Probability and random processes}, Oxford university press, 2020.

\bibitem{hare2020calculus}
{\sc W.~Hare and G.~Jarry-Bolduc}, {\em Calculus identities for generalized simplex gradients: rules and applications}, SIAM Journal on Optimization, 30 (2020), pp.~853--884.

\bibitem{hough2022model}
{\sc M.~Hough and L.~Roberts}, {\em Model-based derivative-free methods for convex-constrained optimization}, SIAM Journal on Optimization, 32 (2022), pp.~2552--2579.

\bibitem{kozak2021stochastic}
{\sc D.~Kozak, S.~Becker, A.~Doostan, and L.~Tenorio}, {\em A stochastic subspace approach to gradient-free optimization in high dimensions}, Computational Optimization and Applications, 79 (2021), pp.~339--368.

\bibitem{larson2019derivative}
{\sc J.~Larson, M.~Menickelly, and S.~M. Wild}, {\em Derivative-free optimization methods}, Acta Numerica, 28 (2019), pp.~287--404.

\bibitem{lewis1999pattern}
{\sc R.~M. Lewis and V.~Torczon}, {\em Pattern search algorithms for bound constrained minimization}, SIAM Journal on optimization, 9 (1999), pp.~1082--1099.

\bibitem{lewis2000pattern}
\leavevmode\vrule height 2pt depth -1.6pt width 23pt, {\em Pattern search methods for linearly constrained minimization}, SIAM Journal on Optimization, 10 (2000), pp.~917--941.

\bibitem{lewis2002globally}
\leavevmode\vrule height 2pt depth -1.6pt width 23pt, {\em A globally convergent augmented lagrangian pattern search algorithm for optimization with general constraints and simple bounds}, SIAM Journal on Optimization, 12 (2002), pp.~1075--1089.

\bibitem{li2021survey}
{\sc Z.~Li, F.~Liu, W.~Yang, S.~Peng, and J.~Zhou}, {\em A survey of convolutional neural networks: analysis, applications, and prospects}, IEEE Transactions on Neural Networks and Learning Systems, 33 (2021), pp.~6999--7019.

\bibitem{meckes2019random}
{\sc E.~S. Meckes}, {\em The random matrix theory of the classical compact groups}, vol.~218, Cambridge University Press, 2019.

\bibitem{menickelly2023avoiding}
{\sc M.~Menickelly}, {\em Avoiding geometry improvement in derivative-free model-based methods via randomization}, 2023.
\newblock \href{https://arxiv.org/abs/2305.17336}{arXiv:2305.17336}.

\bibitem{more1983computing}
{\sc J.~J. Mor{\'e} and D.~C. Sorensen}, {\em Computing a trust region step}, SIAM Journal on Scientific and Statistical Computing, 4 (1983), pp.~553--572.

\bibitem{muirhead2009aspects}
{\sc R.~J. Muirhead}, {\em Aspects of multivariate statistical theory}, John Wiley \& Sons, 2009.

\bibitem{penrose1955generalized}
{\sc R.~Penrose}, {\em A generalized inverse for matrices}, in Mathematical Proceedings of the Cambridge Philosophical Society, vol.~51(3), Cambridge University Press, 1955, pp.~406--413.

\bibitem{Powell1994}
{\sc M.~J.~D. Powell}, {\em A direct search optimization method that models the objective and constraint functions by linear interpolation}, Springer, Dordrecht, 1994, pp.~51--67.

\bibitem{powell2009bobyqa}
\leavevmode\vrule height 2pt depth -1.6pt width 23pt, {\em The {BOBYQA} algorithm for bound constrained optimization without derivatives}, tech. rep., University of Cambridge, 2009.

\bibitem{powell2015fast}
\leavevmode\vrule height 2pt depth -1.6pt width 23pt, {\em On fast trust region methods for quadratic models with linear constraints}, Mathematical Programming Computation, 7 (2015), pp.~237--267.

\bibitem{ragonneau2024pdfo}
{\sc T.~M. Ragonneau and Z.~Zhang}, {\em {PDFO}: a cross-platform package for {P}owell’s derivative-free optimization solvers}, Mathematical Programming Computation, 16 (2024), pp.~535--559.

\bibitem{roberts2025model}
{\sc L.~Roberts}, {\em Model construction for convex-constrained derivative-free optimization}, SIAM Journal on Optimization, 35 (2025), pp.~622--650.

\bibitem{roberts2023direct}
{\sc L.~Roberts and C.~W. Royer}, {\em Direct search based on probabilistic descent in reduced spaces}, SIAM Journal on Optimization, 33 (2023), pp.~3057--3082.

\bibitem{scheinberg2026complexity}
{\sc K.~Scheinberg and A.~Chaudhry}, {\em On complexity of model-based derivative-free methods}, in Proceedings of the International Congress of Mathematicians 2026, vol.~7, SIAM, 2026, pp.~208--228.

\bibitem{shao2021random}
{\sc Z.~Shao}, {\em On random embeddings and their application to optimisation}, PhD thesis, University of Oxford, 2021.

\bibitem{Vershynin_2026}
{\sc R.~Vershynin}, {\em High-Dimensional Probability: An Introduction with Applications in Data Science}, Cambridge Series in Statistical and Probabilistic Mathematics, Cambridge University Press, Cambridge, 2~ed., 2026.

\end{thebibliography}
\end{document}